\documentclass[10pt]{amsart}

\usepackage{epsfig,color}
\usepackage{amsthm}
\usepackage{amssymb}
\usepackage{amsmath,mathrsfs}
\usepackage{epic}
\usepackage{eepic}
\usepackage{graphicx}
\usepackage{psfrag}

\newcommand     {\comment}[1]   {}
\newcommand{\mute}[2] {}
\newcommand     {\printname}[1] {}
\DeclareMathOperator{\sgn}{sgn}%

\newtheorem {Theorem}   {Theorem}
\numberwithin{Theorem}{section}
\newtheorem {Lemma}[equation]    {Lemma}

\newtheorem {Proposition}[equation]{Proposition}
\theoremstyle{definition}
\newtheorem{Definition}[equation]{Definition}
\theoremstyle{remark}
\newtheorem{Remark}[equation]{Remark}
\newtheorem{Example}[equation]{Example}

\def    \eb     {\vec{\mathbf{e}}}

\def    \C      {\mathbb{C}}
\def    \R      {\mathbb{R}}
\def    \Q      {\mathbb{Q}}
\def    \Z      {\mathbb{Z}}

\def    \M      {\mathcal{M}}
\def    \so     {\mathfrak{so}}

\def    \bigds  {\Big/\!\!\Big/}

\newcommand{\dsum}{\displaystyle\sum}  
\newcommand{\dint}{\displaystyle\int}  
\newcommand{\dprod}{\displaystyle\prod}

\begin{document}

\bibliographystyle{plain}

\title[Intersection numbers of polygon spaces]
{Intersection numbers of polygon spaces}

\author[J.\ Agapito]{Jos\'e Agapito}
\author[L.\ Godinho]{Leonor Godinho}
\address{DEPARTAMENTO DE MATEM\'ATICA, INSTITUTO SUPERIOR T\'ECNICO, AV. RO\-VIS\-CO PAIS,
1049-001 LISBON, PORTUGAL, FAX: (351) 21 841 7035}
\email{agapito@math.ist.utl.pt}
\email{lgodin@math.ist.utl.pt}
\thanks{The first author was  partially supported by FCT (Portugal)
through program POCTI/FEDER and grant POCTI/SFRH/BPD/20002/2004; the
second author was partially supported by FCT through program
POCTI/FEDER and  grant POCTI/MAT/57888/2004, and by Funda\c{c}\~{a}o
Calouste Gulbenkian.}

\begin{abstract} We study the intersection ring of the space $\M(\alpha_1,\ldots,\alpha_m)$ of polygons in $\R^3$. We find homology cycles dual to generators of this ring and prove a recursion relation in $m$ (the number of steps)  for their intersection numbers. This result is analog of the recursion relation appearing in the work of Witten and Kontsevich on moduli spaces of punctured curves and on the work of Weitsman on moduli spaces of flat connections on two-manifolds of genus $g$ with $m$ marked points. Based on this recursion formula we obtain an explicit expression for the  computation of  the intersection numbers of polygon spaces and use it in several examples. Among others, we study the special case of equilateral polygon spaces (where all the $\alpha_i$ are the same) and compare our results with the expressions for these particular spaces that have been determined by  Kamiyama and Tezuka. Finally, we relate our explicit formula for the intersection numbers with the generating function for intersection pairings of the moduli space of flat connections of Yoshida, as well as with equivalent expressions for polygon spaces obtained by Takakura and Konno through different techniques. 
\end{abstract}

\maketitle

\section{Introduction}

A Polygon space $\M(\alpha):=\M(\alpha_1,\ldots,\alpha_m)$, $\alpha_i\in \R_+$, is the set of all configurations of closed piecewise linear paths in $\R^3$ with $m$ steps of lengths $\alpha_1,\ldots,\alpha_m$, modulo rotations and translations. These spaces have been widely studied in recent years. For example, Hausmann and Knutson computed their integer cohomology rings in the (generic) case that they are smooth \cite{HK2}. Previously Klyachko \cite{Kl} had showed that the cohomology groups were torsion free and calculated their rank. Moreover, Brion \cite{B} and Kirwan \cite{K} computed the rational cohomology ring for the particular case of equilateral polygons (that is, when the $\alpha_i$'s are all equal) with an odd number of edges as well as this equilateral case modulo the symmetric group. This quotient of the equal-length space by the symmetric group is particularly interesting since it is a compactification of the moduli space of $m$-times-punctured Riemann spheres as well as a compactification of the moduli space of $m$ unordered points in $\C P^1$.

Let us state our main results. In Section~\ref{se:3} we associate a circle bundle $V_i$ on $\M(\alpha)$ to each edge $i=1,\ldots,m$ obtaining a list of degree-$2$ classes: their first Chern classes $c_i:=c_1(V_i)$. It is shown in \cite{HK2}  that these classes generate $H^*(\M(\alpha);\Z[\frac{1}{2}])$. The purpose of our work is the study of their intersection numbers obtaining a recursion relation in $m$ as well as an explicit formula for their computation, and examining its relation with other existing formulas by Konno \cite{Konno}, Takakura \cite{T}, Kamiyama and Tezuka  \cite{KT} and by Yoshida \cite{Y} (this one in the context of moduli spaces of flat $SU(2)$-connections on the $m$-punctured sphere).

As we will see in Section~\ref{se:4} there exist natural homology cycles dual to these Chern classes which will allow us to perform computations in the cohomology ring via intersection theory. In particular, we develop a recursion relation in $m$, the number of edges. For that we start by showing that, for $\alpha_j \neq \alpha_i$, the homology class dual to $c_i$ is given by
$$
[D_{i,j}^+(\alpha)] + \sgn (\alpha_i - \alpha_j) [D_{i,j}^-(\alpha)],
$$ 
where $[D_{i,j}^+(\alpha)]$ and $[D_{i,j}^-(\alpha)]$ are codimension-$2$ submanifolds of $\M(\alpha)$ symplectomorphic to polygon spaces $\M(\alpha_{i,j}^+)$ and $\M(\alpha_{i,j}^-)$ with 
\begin{align*}
\alpha_{i,j}^+ & :=(\alpha_1,\ldots, \hat{\alpha}_i, \ldots, \hat{\alpha}_j,\ldots,\alpha_i + \alpha_j)\\  \alpha_{i,j}^- & :=(\alpha_1,\ldots, \hat{\alpha}_i, \ldots, \hat{\alpha}_j,\ldots,\vert \alpha_i - \alpha_j\vert).
\end{align*} 
These spaces   $\M(\alpha_{i,j}^\pm)$ are themselves endowed with circle bundles and their Chern classes which we denote by $c_i^\pm$. Considering for simplicity $i=m-1$ and $j=m$ (given a permutation $\sigma \in Sym_m$ there is an isomorphism between $\M(\alpha)$ and $\M( \alpha^\sigma)$ given by reordering the steps) we study the behavior of the different classes $c_i$ when restricted to   $[D_{m-1,m}^+(\alpha)]$ and $[D_{m-1,m}^-(\alpha)]$ and compare it with the behavior of $c_i^\pm$. In particular, we show that
\begin{align*}
& (i_\pm \circ s_\pm^{-1})^* c_i    =  c_i^\pm \quad \text{for}\quad 1 \leq i\leq m-2; \\
&(i_+ \circ s_+^{-1})^* c_{m-1}    =  c_{m-1}^+; \\
&(i_- \circ s_-^{-1})^* c_{m-1}    =  \sgn (\alpha_{m-1}-\alpha_m) \,c_{m-1}^-; \\
&(i_+ \circ s_+^{-1})^* c_{m}     =  c_{m-1}^+ ;   \\
& (i_- \circ s_-^{-1})^* c_{m}    =  - \sgn (\alpha_{m-1}-\alpha_m)\,  c_{m-1}^-,
\end{align*} 
where $i_{\pm}:D_{m-1,m}^\pm(\alpha)  \hookrightarrow \M(\alpha)$ denote the inclusion maps and $s_{\pm}:D_{m-1,m}^\pm(\alpha) \to  \M(\alpha_{i,j}^\pm)$ are the symplectomorphisms which identify these spaces. We then obtain the following recursion relation (cf. Section~\ref{se:5}).
\begin{Theorem}\label{thm:1.1}
Suppose $\alpha_m\neq \alpha_{m-1}$ and let  $c_i^+$ and  $c_i^-$  be the cohomology classes $c_1(V_m(\alpha^+))$ and $c_1(V_m(\alpha^-))$, where
$$
\alpha^+:=(\alpha_1,\ldots,\alpha_{m-2},\alpha_{m-1} +  \alpha_m) \quad \text{and} \quad \alpha^-:=(\alpha_1,\ldots,\alpha_{m-2},\vert \alpha_{m-1} - \alpha_m\vert ).
$$
Then, for $k_1,\ldots, k_m \in \Z_{\geq 0}$ such that $k_1+\cdots + k_m= m-3$ and $k_m\geq 1$,
\begin{equation}
\begin{split}
&\dint\limits_{\M(\alpha)} c_1^{k_1} \cdots c_m^{k_m} = 
\dint\limits_{\M(\alpha^+)}
 (c^+_1)^{k_1} \cdots (c^+_{m-2})^{k_{m-2}} (c_{m-1}^+)^{k_{m-1}+k_m - 1} \quad + \\
  \\
& + \,  (-1)^{k_m - 1} \left( \sgn( \alpha_{m-1} -\alpha_m ) \right)^{k_{m-1}+k_m}
 \dint\limits_{\M(\alpha^-)}\!\! (c^-_1)^{k_1}\cdots
(c^-_{m-2})^{k_{m-2}} (c^-_{m-1})^{k_{m-1}+k_m - 1}. 
\end{split}
\end{equation}
\end{Theorem}  
This recursion relation is analog of the recursion relation appearing in the work of Witten and Kontsevich on moduli spaces of punctured curves \cite{Kont,Wi1,Wi3} and on the work of Weitsman on moduli spaces of flat connections on two-manifolds of genus $g$ with $m$ marked points \cite{W}.  This is not surprising since for small values of $\alpha$ the polygon spaces $\M(\alpha)$ can be identified with moduli spaces of flat $SU(2)$-connections on the $m$-punctured sphere (cf. Section~\ref{se:9}). The proof of Theorem~\ref{thm:1.1} takes profit from this identification as it follows Weitsman's proof \cite{W} of the corresponding recursion relation for moduli spaces of flat connections on surfaces of genus $g$, adapting it to the context of polygon spaces and making the necessary changes for the genus $g=0$ situation. 

Based on this recursion relation we obtain an explicit expression for the  computation of  the intersection numbers (cf. Section~\ref{se:6}).
\begin{Theorem}\label{thm:1.2}
Let $\alpha=(\alpha_1,\ldots, \alpha_m)$ be generic. Suppose $k_{m-l},\ldots, k_m\in \Z_+$,  $k_1=\cdots = k_{m-l-1}=0$ and $k_{m-l} + \cdots + k_m = m-3$. Let  $c_i:=c_1(V_i(\alpha))$ be the first Chern classes of the circle bundles $V_i(\alpha)\to \M(\alpha)$. Then
\begin{equation}\nonumber
\int_{\M(\alpha)} c_{m-l}^{k_{m-l}} \cdots c_m^{k_m} = \sum_{J \in \mathcal{T}(\alpha)} (-1)^{\left( \sum_{i \in I\setminus J} k_i \right) + m - \vert J  \vert}.
\end{equation}
\end{Theorem}
Here $\mathcal{T}(\alpha)$ is a special family of sets $J\subset I:=\{3,\ldots,m\}$ that we will call \emph{triangular} for which $
\sum_{j=3}^m (-1)^{\chi_{I\setminus J}(i)} \, \alpha_i > 0 $
and  for which the following triangle inequalities are satisfied:
\begin{align*}
&\alpha_1  \leq \alpha_2 + \sum_{i=3}^m (-1)^{\chi_{I\setminus J}(i)} \, \alpha_i, \\
&\alpha_2  \leq \alpha_1 + \sum_{i=3}^m (-1)^{\chi_{I\setminus J}(i)} \, \alpha_i, \\
& \sum_{i=3}^m (-1)^{\chi_{I\setminus J}(i)} \, \alpha_i  \leq \alpha_1 + \alpha_2.
\end{align*}

After we prove this relation we work out several examples in Section~\ref{se:7} and compare the results obtained through our formulas with the ones obtained by Hausmann and Knutson in \cite{HK2}.
Then, in Section~\ref{se:8}, we study the equilateral case $\M_m$. Requiring $\alpha$ to be generic in this case means exactly that $m$ is odd. The intersection numbers for these particular spaces  have  been determined by Kamiyama and Tezuka \cite{KT}. 
\begin{Theorem}(Kamiyama-Tezuka)
Let $(d_1,\ldots,d_m)$ be a sequence of nonnegative integers with $\sum d_i=m-3$. Let $\beta_i, \varepsilon_i$ be such that $d_i=2\beta_i + \varepsilon_i$, where $\varepsilon_i=0$ or $1$. Then, defining
$$
\rho_{m,2k}:= (-1)^k \frac{\left( \begin{array}{c} \frac{m-3}{2} \\ k \end{array} \right) \left(\begin{array}{c} m-2 \\ \frac{m-1}{2} \end{array}\right)} { \left( \begin{array}{c} m-2  \\  2k + 1 \end{array} \right)},
$$
we have
\begin{enumerate}
\item if $\beta_i=0$ for $1\leq i \leq m$ then 
$\int c_1^{d_1}\cdots c_m^{d_m} = \rho_{m,0}$;
\item[]
\item if $\beta_i \neq 0$ for some $i$ then 
$\int c_1^{d_1}\cdots c_m^{d_m} = \rho_{m,2k}$, with $k=\beta_1+\cdots+ \beta_m$.
\end{enumerate}
\end{Theorem}
We use Theorem~\ref{thm:1.2} to compute these numbers and we show the equivalence between our  results and the ones obtained by Kamiyama and Tezuka, through several combinatorial computations. For completion,  we consider the quotient of $\M_m$ by the action of the symmetric group. The cohomology ring of this space was computed by Brion \cite{B}, Klyachko \cite{Kl} and Hausmann and Knutson  \cite{HK2}. It can be identified with the invariant part of $H^*(\M_m,;\Q)$ by the symmetric group action and is generated by $\sigma_1$, the first invariant symmetric polynomial in the classes $c_i$ and, for example, by $c_m^2$, a generator of degree $4$ (or any other $c_i^2$ since they are all equal to the Pontrjagin class of the principal $SO(3)$-bundle $A_m \to \M_m$, where $A_m:=\{(\vec{v}_1,\ldots,\vec{v}_m)\in (\R^3)^m\mid \sum_{i=1}^m \vec{v}_i = 0,\,\, \vert \vec{v}_1\vert=\cdots = \vert \vec{v}_m\vert\}$ ). We compute the intersection numbers $\sigma_1^k \cdot c_m^{m-3-k}$ for an even $k$ and obtain the following result.
\begin{Proposition} For an even integer $k$
\begin{align*}
& \int_{\M_m/Sym_m}  \sigma_1^k \cdot c_m^{m-3-k}  = \\ &=  (-1)^{\frac{m-3}{2}}  \left( \begin{array}{c} m-2 \\ \frac{m-1}{2} \end{array}\right)  \sum_{j=0}^{\frac{k}{2}} (-1)^j \frac{ \left( \begin{array}{c} \frac{m-3}{2} \\ j  \end{array} \right)}{ \left( \begin{array}{c} m-2 \\ 2j \end{array}\right)} \sum_{ \small{\begin{array}{c} k_1+\cdots+k_m=k \, \text{s.t.} \\ 2j\, \text{of the $k_i$'s are odd}\end{array}}} \left(\begin{array}{c} k \\ k_1, \cdots, k_m \end{array} \right).
\end{align*}
\end{Proposition}

In Section~\ref{se:9}, we explore the identification between polygon spaces and moduli spaces of flat $SU(2)$-connections on the $m$-punctured sphere and we relate the generating function for intersection pairings obtained by Yoshida \cite{Y} in the context of the moduli space of flat connections with our explicit formula from Theorem~\ref{thm:1.2}.

Finally, in Section~\ref{se:10} we compare Theorem~\ref{thm:1.2} with the formula for intersection pairings on polygon spaces previously obtained by Takakura \cite{T} and Konno \cite{Konno} by different methods. In particular, Takakura in \cite{T} uses the ``quantization commutes with reduction'' theorem of Guillemin and Sternberg \cite{GS} to obtain an explicit formula for intersection pairings on polygon spaces, and Konno in \cite{Konno} computes these pairings from an algebro-geometric point of view generalizing the methods of Kamiyama and Tezuka \cite{KT}. The two expressions are equivalent (see \cite{Konno}) although Konno uses a different basis, 
$$
v_i := \frac{c_i + c_m}{2},\, i=1,\ldots, m
$$
to write his formula. Their result is the following.
\begin{Theorem}(Konno, Takakura)

Let $\alpha=(\alpha_1,\ldots, \alpha_m)$ be generic  and let
$(k_1,\ldots,k_m)$ be a sequence of nonnegative integers with $\sum_{i=1}^m k_i=m-3$. Let $\mathcal{S}(\alpha)$ be the family of sets  $R \subset \{1, \ldots, m \}$ for which $\sum_{i \in R} \alpha_i - \sum_{i \notin R} \alpha_i < 0$. Then we have
$$
\int_{\M(\alpha)} c_1^{k_1} \cdots c_m^{k_m} = -\frac{1}{2} \sum_{R \in \mathcal{S}(\alpha)} (-1)^{\vert R \vert + \sum_{i \in R} k_i}.
$$
\end{Theorem}
We will show that this formula is equivalent to ours. Note, however, that our formula uses a smaller  family of sets (the  triangular sets $\mathcal{T}(\alpha)$) which is contained in $\mathcal{S}(\alpha)$. 

Moreover, Konno's method (like the one used in \cite{KT}) relies on the explicit relations among the generators of the cohomology ring. Instead, our approach and Takakura's do not use these relations and contain enough information to recover the structure of this ring.

\noindent{\bfseries Acknowledgements.} The authors are grateful to Jonathan Weitsman for suggesting this problem and for many useful conversations and remarks. They would also like to thank Jean-Claude Hausmann for his comments on an earlier version of this work and to Takakura and Konno for bringing their work to our attention. 
\bigskip

\section{Setting and a brief survey on polygon spaces}
Let $\alpha=(\alpha_1,\ldots,\alpha_m)\in\R^m_+$ and, for each $1\le
i\le m$, let $S^2_{\alpha_i}$ be the sphere of radius $\alpha_i$ in $\R^3$. The product $\prod_{i=1}^m
S^2_{\alpha_i}$ can be thought of as the space of all paths starting
at the origin with $m$ consecutive steps $\vec{v}_i$ whose lengths
are $\vert\vec{v}_i\vert=\alpha_i$. We are interested in closed polygonal paths, i.e.,
$$
\left\{\vec{v}=(\vec{v}_1,\ldots,\vec{v}_m)\in\prod_{i=1}^m
S^2_{\alpha_i}\mid \,\, \sum_{i=1}^m \vec{v}_i = 0 \right\}.
$$
The polygon
space $\M(\alpha)$ is the space of all configurations of closed polygonal
paths in $\R^3$ starting and ending at the origin with $m$ edges of lengths $\alpha_1,\ldots,\alpha_m$,
modulo rotation. More precisely,
$$
\M(\alpha):=\left\{\vec{v}=
(\vec{v}_1,\ldots,\vec{v}_m)\in\prod_{i=1}^m
S^2_{\alpha_i}\mid \,\, \sum_{i=1}^m \vec{v}_i = 0 \right\}\Big/
SO(3).
$$
We want to rule out \textit{degenerate polygons}, that is, polygons which are
contained in a straight line. This degeneracy condition is equivalent to finding
$\varepsilon_i=\pm 1$ for $1\le i\le m$, such that
\begin{equation}
\label{eq:degeneracy}
\sum_{i=1}^m \varepsilon_i\alpha_i=0.
\end{equation}
Hence, we say that $\alpha$ is \textit{generic} if equation
\eqref{eq:degeneracy} has no solution with $\varepsilon_i=\pm 1$. If
$\alpha$ is generic, the quotient space $\M(\alpha)$ is a symplectic
manifold of dimension $2(m-3)$. Indeed, considering the product
symplectic structure
$\omega=\alpha_1\omega_{S^2}+\ldots+\alpha_m\omega_{S^2}$ on
$\prod_{i=1}^m S^2_{\alpha_i}\subset(\R^3)^m$, where $\omega_{S^2}$
is the standard symplectic structure on $S^2$, and the diagonal
Hamiltonian $SO(3)$-action with moment map
$$
\begin{array}{rccl}
\mu\colon & \dprod_{i=1}^m S^2_{\alpha_i} & \longrightarrow &
\so(3)^*\cong(\R^3)^* \\[2ex]
          & \vec{v}=(\vec{v}_1,\ldots,\vec{v}_m) & \longmapsto &
          \mu(\vec{v})=\dsum_{i=1}^m \vec{v}_i\quad \mbox{(``endpoint'')},
\end{array}
$$
we can, for a generic $\alpha$, see the polygon space $\M(\alpha)$
as a symplectic quotient of the path space. When $\alpha$ is generic the $SO(3)$-action is free and $0$ is a regular value of $\mu$ and so $\M(\alpha)$ is itself a symplectic
manifold. More precisely,
$$
\M(\alpha) = \mu^{-1}(0) \Big/ SO(3) = \left(\prod_{i=1}^m
S^2_{\alpha_i}\right) \underset{0\,\,\,}{\bigds}\,
SO(3)\cong\prod_{i=1}^{m-1}S^2_{\alpha_i}\underset{\alpha_m\,\,\,}{\bigds}\,
SO(3),
$$
where this last space is the set of paths of $m-1$ steps of lengths
$\alpha_1\,\ldots,\alpha_{m-1}$ whose endpoint is at a distance
$\alpha_m$ from the origin modulo $SO(3)$ (see \cite{Kl,KM,HK2} for additional details).  Hereafter, we will assume
that $\alpha$ is generic.

\begin{Remark}
\label{rk:2.1}
For $\alpha$ and $\alpha^\prime$ generic and sufficiently close, $\M(\alpha)$ and $\M(\alpha^\prime)$ are diffeomorphic (\cite{HK2}, Proposition 2.2). Indeed, if they are close enough such that $\alpha(t):= t\alpha^\prime + (1-t)\alpha$ is always generic for $t\in [0,1]$, the map 
\begin{align*}
\beta: [0,1]\times \prod_{i=1}^m S^2 & \to [0,1] \times \R^3 \\
(t, z_1,\ldots,z_m) & \mapsto (t, \mu(\alpha_1(t)z_1,\ldots, \alpha_m(t)z_m))
\end{align*}
has no critical values in $[0,1]\times \{0\}$ and so its gradient flow gives an $SO(3)$-equivariant diffeomorphism between
$\beta^{-1}(0,0)$ and $\beta^{-1}(1,0)$. Hence, the quotient spaces $\M(\alpha)=\beta^{-1}(0,0)/SO(3)$ and   $\M(\alpha^\prime)=\beta^{-1}(1,0)/SO(3)$ are also diffeomorphic.
\end{Remark}
The restriction of the diagonal $SO(3)$-action on the path space to an action of $S^1\cong SO(2)$ is also Hamiltonian with moment map 
$$
\begin{array}{rccl}
\overline{\mu}\colon & \dprod_{i=1}^m S^2_{\alpha_i} & \longrightarrow &
\so(2)^* \cong(\R)^* \\[2ex]
          & \vec{v}=(\vec{v}_1,\ldots,\vec{v}_m) & \longmapsto &
          \overline{\mu}(\vec{v})=\zeta\left(\dsum_{i=1}^m \vec{v}_i \right) \quad \mbox{(``height of endpoint'')},
\end{array}
$$
where $\zeta$ is the projection $\zeta(x,y,z)=z$.

Using the maps $\zeta$ and $\overline{\mu}$ we will define some other spaces that will be relevant to us. 
The first is the \emph{abelian polygon space},  
$$
\mathcal{AM}(\alpha)  := \left\{\vec{v}=
(\vec{v}_1,\ldots,\vec{v}_{m-1})\in\prod_{i=1}^{m-1}
S^2_{\alpha_i}\mid \,\, \zeta(\sum_{i=1}^{m-1} \vec{v}_i) = \alpha_m \right\}\Big/
S^1,
$$
the space of piecewise linear $(m-1)$-chains with edge lengths $\alpha_1,\ldots, \alpha_{m-1}$ which end on the plane $z=\alpha_m$, modulo rotations around the $z$-axis. It is a symplectic quotient of the path space  $\prod_{i=1}^{m-1}
S^2_{\alpha_i}$ by the maximal torus $S^1$ of $SO(3)$,
\begin{equation}
\label{eq:2.1}
\mathcal{AM}(\alpha)=\zeta^{-1}(\alpha_m)/S^1,
\end{equation}
and so it is a symplectic manifold of dimension $2(m-2)$ containing $\mathcal{M}(\alpha)$ as a symplectic submanifold of codimension $2$. 

\begin{Remark}
\label{rk:2.2}
Note that we can always rotate any element $\vec{v}\in \M(\alpha)$ in such a way that the sum $\sum_{i=1}^{m-1} \vec{v}_i$ ends not only on the plane $z=\alpha_m$ but also on the $z$-axis, so that $\vec{v}_m$ points downwards. 
\end{Remark}

The \emph{upper path space} $\mathcal{UP}(\alpha)$ is defined as the space
$$
\mathcal{UP}(\alpha):= 
\left\{ \vec{v} = (\vec{v}_1,\ldots, \vec{v}_{m-1}) \in \prod_{i=1}^{m-1} S^2_{\alpha_i} \mid \,\, \zeta \left(\sum_{i=1}^{m-1} \vec{v}_i\right) \geq \alpha_m \right\}\Big/ \sim,
$$
where $\vec{v} \sim \vec{v}^\prime$ if and only if $\vec{v}= \vec{v}^\prime$ or $ \zeta(\sum_{i=1}^{m-1} \vec{v}_i)= \alpha_m$ and $[\vec{v}]=[\vec{v}^\prime]$ in $\mathcal{AM}(\alpha)$. It is the symplectic cut (in the sense of \cite{L}) of the path space $\prod_{i=1}^{m-1} S^2_{\alpha_i}$ at the level $\alpha_m$ of the moment map $\overline{\mu}$. Hence, it is a symplectic manifold of dimension $2(m-1)$ which  contains $\mathcal{AM}(\alpha)$ as a submanifold of codimension $2$.

\section{Circle bundles over polygon spaces}\label{se:3}

We construct circle bundles over $\M(\alpha)$ as follows. We define, for
each $1\le j\le m$, the space
$$
V_{j}(\alpha):=\left\{\vec{v}=(\vec{v}_1,\ldots,\vec{v}_m)\in\prod_{i=1}^m
S^2_{\alpha_i}\mid \,\, \sum_{i=1}^m \vec{v}_i = 0
\,\,\mbox{and}\,\, \vec{v}_j=(0,0,\alpha_j)\right\},
$$
which is a smooth manifold of dimension $2m-5$. The circle $S^1$
acts on $V_{j}(\alpha)$ by rotation around the $z$-axis. As $\alpha$
is generic, this action is free and it is easy to check that
$V_j(\alpha) \Big/ S^1=\M(\alpha)$. Hence, $V_j(\alpha)\to\M(\alpha)$
is a principal circle bundle determined by its Chern class 
$$
c_j:=c_1(V_j(\alpha))\in H^2(\M(\alpha); \Z).
$$ 

Let us also consider the level set $B(\alpha):=  \zeta^{-1}(\alpha_m)\subset \prod_{i=1}^{m-1}
S^2_{\alpha_i}$, i.e.
$$
B(\alpha):=\left\{\vec{v}=(\vec{v}_1,\ldots,\vec{v}_{m-1})\in \prod_{i=1}^{m-1}
S^2_{\alpha_i} \mid \,\, \zeta\left(\sum_{i=1}^{m-1} \vec{v}_i\right) = \alpha_m \right\}.
$$
Since $\alpha$ is generic, $S^1$ acts freely on this space making it a principal bundle over $\mathcal{AM}(\alpha)$.
Moreover, we have the following commutative diagram 
$$
\begin{array}{ccc}
V_m(\alpha) & \stackrel{\tilde{i}}{\longrightarrow} & B(\alpha) \\ \downarrow & & \downarrow\\
\M(\alpha) & \stackrel{i}{\longrightarrow}  & \mathcal{AM}(\alpha). \\
\end{array}
$$
Note that the inclusion $\tilde{i}:V_m(\alpha) \hookrightarrow B(\alpha)$ is anti-equivariant (i.e. $\tilde{i}(\lambda \cdot \vec{v})=\lambda^{-1}\cdot  \tilde{i}(\vec{v})$) since, in the identification of $\M(\alpha)$ as a submanifold of $\mathcal{AM}(\alpha)$, the vector $\vec{v}_m$ must face downward (see Remark~\ref{rk:2.2}). Hence, $c_m=-i^*(c_1(B(\alpha)))$.

Let $\mathcal{N}$ be a tubular neighborhood of $\mathcal{AM}(\alpha)$ inside $\mathcal{UP}(\alpha)$. The retraction $\mathcal{N} \to \mathcal{AM}(\alpha)$ is the disc bundle associated to the circle bundle $B(\alpha) \to \mathcal{AM}(\alpha)$ and so $c_1(B(\alpha))$ is the Euler class of the normal bundle of  $\mathcal{AM}(\alpha)$ inside  $\mathcal{UP}(\alpha)$. On the other hand, since $\mathcal{AM}(\alpha)$ is the reduced space \eqref{eq:2.1},
$$
\zeta^{-1}(\alpha_m) \Big/ S^1,
$$
we have, by the Duistermaat-Heckmann theorem \cite{DH,G}, that 
$$
c_1(B(\alpha))=- \frac{\partial}{\partial \alpha_m} [\omega]
$$
in $H^2(\mathcal{AM}(\alpha);\R)$, and so 
$$
c_m=\frac{\partial}{\partial \alpha_m} [\omega]
$$ 
in  $H^2(\mathcal{M}(\alpha); \R)$. By symmetry, since any edge can be the last, we have 
\begin{equation}\label{eq:3.1}
c_j=\frac{\partial}{\partial \alpha_j} [\omega].
\end{equation}
It is shown in \cite{HK2} that these classes $c_j$ generate $H^*(\M(\alpha); \Z[\frac{1}{2}])$.

\begin{Remark}
\label{rk:3.1}
Given a permutation $\sigma\in \text{Sym}_m$ there is an isomorphism between $M(\alpha)$ and $\M(\alpha^\sigma)$ given by reordering the steps (note that a polygon is simply a list of $m$ vectors whose sum is zero modulo rotation). From the above geometric construction of the bundles $V_j(\alpha)$ we see that the induced isomorphism on cohomology $H^2(\M(\alpha^\sigma)) \to H^2(\M(\alpha))$ yields $c_j\to c_{\sigma(j)}$. 
\end{Remark}

\section{Dual homology classes}\label{se:4}
The content of the next two sections takes profit from the identification, for small values of $\alpha$, between polygon spaces $\M(\alpha)$ and moduli spaces of flat $SU(2)$- connections on the $m$-punctured sphere (cf. Section~\ref{se:9}). Indeed,  we follow the proof of Weitsman's recursion relation \cite{W}  for moduli spaces of flat connections on surfaces of genus $g$, adapting it to the context of polygon spaces and making the necessary changes for the genus $g=0$ situation.  

We are interested in determining homology cycles representing the
Poincar\'e duals of the first Chern classes
$c_j:=c_1(V_j(\alpha)) \in H^2(\M(\alpha); \Z)$. For this purpose,
let $i$ and $j$ be positive integers such that $1\le i,j\le m$ and
$i\neq j$. Let us denote by $D_{i,j}(\alpha)$ the submanifold of
$\M(\alpha)$ formed by those polygons $\vec{v}\in\M(\alpha)$ for which
the step $\vec{v}_i$ is parallel to $\vec{v}_j$. By rotating
$\vec{v}$ around a suitable axis, if necessary, we can assume that
$\vec{v}_i$ and $\vec{v}_j$ are parallel to the $z$ axis. Hence
$$
D_{i,j}(\alpha):=\left\{\vec{v}=(\vec{v}_1,\ldots,\vec{v}_m)\in\M(\alpha)
  \mid \,\, \vec{v}_i=\alpha_i \eb_3, \,\,  \vec{v}_j= \alpha_j \eb_3\right\}
$$
is a codimension two submanifold of $\M(\alpha)$. We will show that the bundle $V_j(\alpha)$ has a section on the complement of $D_{i,j}(\alpha)$ and so $c_j$ will be dual to an element of $i_*(H_{2(n-4)}(D_{i,j}(\alpha)))$, where $i:D_{i,j}(\alpha) \hookrightarrow \M(\alpha)$ is the inclusion map. Note that  $D_{i,j}(\alpha)=D_{j,i}(\alpha)$ so  $D_{i,j}(\alpha)$  will not be connected. Indeed, both $V_i(\alpha)$ and  $V_j(\alpha)$ are trivial on the complement of $D_{i,j}(\alpha)$ but they are not powers of the same circle bundle.

\begin{Proposition}
The circle bundle $V_j(\alpha)\vert_{\M(\alpha)\setminus D_{i,j}(\alpha)}\stackrel{\pi_j}{\longrightarrow} \M(\alpha)\setminus D_{i,j}(\alpha)$ has a section.
\end{Proposition}

\proof
A section $s:\M(\alpha)\setminus D_{i,j}(\alpha) \to V_j(\alpha)$ can be given by choosing, for each $\vec{v} \in \M(\alpha)$, a representative $s( \vec{v})\in V_j(\alpha)$ such that $\pi_j(s(\vec{v}))=\vec{v}$. To do this we assign to each $\vec{v}$ the unique element $p\in \pi_j^{-1}(\vec{v}) \subset V_j(\alpha)$ for which the $i$-th step $\vec{v}_i$ projects onto the $xOy$ plane along the positive $x$-axis. Such a representative will always exist in  $\pi_j^{-1}(\vec{v})$ as long as $\vec{v}_i$ is not parallel to $\vec{v}_j$, that is, as long as $\vec{v}\notin D_{i,j}(\alpha)$.
\endproof

Since $D_{i,j}(\alpha)$ consists of classes of polygons for which the step $\vec{v}_i$ is parallel to $\vec{v}_j$, it clearly has  two different connected components 
\begin{align*}
& D_{i,j}^+(\alpha) = \left\{\vec{v}\in D_{i,j}(\alpha) \mid 
\langle\vec{v}_i,\vec{v}_j \rangle > 0\right\} \quad \text{and} \quad D_{i,j}^-(\alpha) = \left\{\vec{v}\in D_{i,j}(\alpha)\mid
\langle\vec{v}_i,\vec{v}_j\rangle < 0\right\}.
\end{align*}
Both spaces $D_{i,j}^\pm(\alpha)$ are symplectomorphic to polygon spaces obtained by permuting the edges so that the steps $\vec{v}_i$ and $\vec{v}_j$ become the last two consecutive edges. In other words,

\begin{Proposition}\label{le:generic}
For $\alpha_i \neq \alpha_j$ there exist symplectomorphisms  
$$
\begin{array}{cccl}
s_\pm: & D_{i,j}^\pm(\alpha)& \longrightarrow & \M(\alpha_{i,j}^\pm)\\
 & [(\vec{v}_1,\ldots,\vec{v}_m)] & \mapsto & [(\vec{v}_1,\ldots, \hat{\vec{v}}_i, \ldots, \hat{\vec{v}}_j, \ldots, \vec{v}_i \pm \vec{v}_j)],
\end{array}
$$
where 
\begin{eqnarray*}
\alpha_{i,j}^+ & = & \left(\alpha_1,\ldots,\hat{\alpha_i}, \ldots, \hat{\alpha_j}, \ldots, \alpha_i + \alpha_j\right),\\
\alpha_{i,j}^- & = & \left(\alpha_1,\ldots,\hat{\alpha_i},\ldots, \hat{\alpha_j},\ldots, \vert\alpha_i-\alpha_j\vert \right).
\end{eqnarray*} 
\end{Proposition}
Note that both $\alpha_{i,j}^+$ and $\alpha_{i,j}^-$ are generic provided that $\alpha$ is. We conclude that the manifolds $D_{i,j}^\pm(\alpha)$ are connected symplectic manifolds and so we can orient them using the symplectic form by requiring that $\int_{D_{i,j}^\pm(\alpha)} \omega^{n-4} > 0$. We obtain in this way two generators of $H_{2(n-4)}(D_{i,j}^\pm(\alpha))$: $[D_{i,j}^+(\alpha)]$ and  $[D_{i,j}^-(\alpha)]$. Therefore, to determine the Poincar\'e dual of $c_j$ we just have to determine two constants $a_{i,j}$ and $b_{i,j}$ as in the following proposition.
\begin{Proposition}
Let $i: D_{i,j}(\alpha)\hookrightarrow \M(\alpha)$ be the inclusion map. If  $\alpha_i \neq \alpha_j$ then the Poincar\'e dual of $c_j$ is in $i_* H_{2(n-4)}(D_{i,j}(\alpha))$ and can be written as
$$
a_{i,j}[D_{i,j}^+(\alpha)] + b_{i,j} [D_{i,j}^-(\alpha)].
$$
\end{Proposition}
We will now see that the constants  $a_{i,j}$ are always equal to $1$ while the constants $b_{i,j}$ depend on the values of $\alpha_i$ and $\alpha_j$.
\begin{Proposition}
\label{prop:4.1}
The  constants  $a_{i,j}$ and $b_{i,j}$ are equal to 
$$
a_{i,j}=1\quad \text{and} \quad b_{i,j}=\sgn(\alpha_i-\alpha_j).
$$
\end{Proposition}

\proof
For simplicity, we consider $i=m-1$ and $j=m$ (we can always reduce to this case due to Remark~\ref{rk:3.1}). Take a fixed polygon $\vec{\underline{v}}$ in $\R^3$ such that $\vert \vec{\underline{v}}_i \vert =\alpha_i$ and consider the space
\begin{align*}
W := \Bigl\{ \vec{v}=(\vec{v}_1,\ldots, \vec{v}_m)\in \Bigr. & (\R^3)^m \mid  \quad \sum_{k=1}^m \vec{v}_k = 0, \quad \vert \vec{v}_k \vert = \alpha_k, \quad k=1,\ldots,m  \\  & \Bigl. \text{and} \quad \vec{v}_i=\vec{\underline{v}}_i, \quad  i=1,\ldots, m-3 \Bigr\}\Big/ SO(3).
\end{align*} 
This corresponds to fixing the first $m-3$ edges of the polygons, allowing to move only the last three. The manifold $W$ is symplectomorphic to the sphere $\M(l,\alpha_{m-2},\alpha_{m-1},\alpha_m)$, where $l:=\vert \sum_{k=1}^{m-3} \vec{\underline{v}}_k \vert$. To find the constants $a_{i,j}$ and $b_{i,j}$ we just have to compare $c_m\vert_W$ with $W\cap D_{m-1,m}^\pm(\alpha)$. 

We compute  $c_m\vert_W$ using the variation of the symplectic form given in \eqref{eq:3.1}. Hence, if $i: W \hookrightarrow \M(\alpha)$ is the inclusion map then 
$$
c_m\vert_W=i^*c_m= i^* \frac{\partial}{\partial \alpha_m} [\omega]= \frac{\partial}{\partial \alpha_m} [i^* \omega].
$$
On the other hand, $W$ is a toric manifold (equivariantly symplectomorphic to a sphere equipped with a circle action) with moment map 
$$\mu:W \cong \M(l,\alpha_{m-2},\alpha_{m-1},\alpha_m) \to \R^*$$ 
obtained by restricting to $W$ the action on $\M(\alpha)$ given by bending the polygons along the last diagonal. The moment map image $\mu(W)$ is the interval
$$
\Delta = \left[\, \max{\{\vert l - \alpha_{m-2}\vert,\, \vert \alpha_{m-1} - \alpha_m\vert\}},\, \min{\{l + \alpha_{m-2},\, \alpha_{m-1} + \alpha_m}\}\, \right]
$$
(see Section $6.$ of  \cite{HK1} for details). 

Now let $M^{2n}$ be an arbitrary  toric manifold with moment map $\mu$, equipped with a family  of symplectic forms $\Omega_t$. The corresponding family of moment polytopes $\Delta_t\subset \frak{t}^*$,  is determined by their facets 
$$
F_k=\{x\in \R^n \mid \langle x, u_k \rangle = \lambda_k(t) \}, \quad k=1,\ldots,N,
$$
where $N$ is the number of facets of $\Delta_t$.
Suppose that the polytopes $\Delta_t$ stay  combinatorially the same as $t$ changes but the value  $\lambda_i(t)$ for some $i\in \{1,\ldots, N \}$ depends linearly on $t$ and, as $t$ increases, the facet $F_i$ moves outward while the others stay fixed. Then we know that $\frac{d}{d t} \Omega_t$ is the Poincar\'e dual of the homology class $[\mu^{-1}(F_i)]$ where the orientation of $\mu^{-1}(F_i)$ is given by requiring $\int_{\mu^{-1}(F_i)} \Omega_t^{n-1} > 0$ (cf. \cite{G} for details).

Applying this general fact about toric manifolds to $W$ and its moment polytope $\Delta$, we see that as $\alpha_m$ changes, the cohomology class of the symplectic form of $W$, $[\omega_{\alpha_m}]$, changes by the Poincar\'e dual of the homology class
\begin{align*}
& [\mu^{-1}(\alpha_{m-1}+\alpha_m)]  + \sgn(\alpha_{m-1}-\alpha_m) [\mu^{-1}(\vert \alpha_{m-1} - \alpha_m\vert)] = \\  & = [D_{m-1,m}^+(\alpha) \cap W] +  \sgn(\alpha_{m-1}-\alpha_m) [D_{m-1,m}^-(\alpha) \cap W]
\end{align*}
and the result follows.
\endproof

\section{A recursion formula}\label{se:5}

In order to prove Theorem~\ref{thm:1.1} we still need to study the behavior of the different Chern classes $c_i$ when restricted to $[D_{m-1,m}^+(\alpha)]$ and  $[D_{m-1,m}^-(\alpha)]$.
\begin{Proposition}
\label{prop:5.1} Suppose $\alpha_m\neq \alpha_{m-1}$ and
let $c_m^+$ and $c_m^-$ be the cohomology classes
$c_1\left(V_m(\alpha^+)\right)$ and $c_1\left(V_m(\alpha^-)\right)$, where 
$$
\alpha^+:=(\alpha_1,\ldots,\alpha_{m-2},\alpha_{m-1} +  \alpha_m) \quad \text{and} \quad \alpha^-:=(\alpha_1,\ldots,\alpha_{m-2},\vert \alpha_{m-1} - \alpha_m\vert ).
$$
Then, considering  the inclusion maps $i_\pm: D_{m-1,m}^\pm(\alpha) \hookrightarrow \M(\alpha)$ and the symplectomorphisms  
$s_\pm:  D_{i,j}^\pm(\alpha) \longrightarrow  \M(\alpha_{i,j}^\pm)$ from Proposition~\ref{le:generic}, we have
\begin{align*}
& (i_\pm \circ s_\pm^{-1})^* c_i    =  c_i^\pm \quad \text{for}\quad 1 \leq i\leq m-2; \\
&(i_+ \circ s_+^{-1})^* c_{m-1}    =  c_{m-1}^+; \\
&(i_- \circ s_-^{-1})^* c_{m-1}    =  \sgn (\alpha_{m-1}-\alpha_m) \,c_{m-1}^-; \\
&(i_+ \circ s_+^{-1})^* c_{m}     =  c_{m-1}^+ ;   \\
& (i_- \circ s_-^{-1})^* c_{m}    =  - \sgn (\alpha_{m-1}-\alpha_m)\,  c_{m-1}^-.  
\end{align*}
\end{Proposition}
\proof 
Let us denote by $\omega^+$ and $\omega^-$ the symplectic forms on $\M(\alpha^+)$ and $\M(\alpha^-)$ obtained by restricting to $D_{m-1,m}^\pm(\alpha)$ the symplectic form $\omega$ of $\M(\alpha)$. Then, for $1 \leq i\leq m-2$, we have
$$
c_i^\pm = \frac{\partial}{\partial \alpha^\pm_i}[\omega^\pm] = \frac{\partial}{\partial \alpha_i}[(i_\pm \circ s_\pm^{-1})^* \omega] = (i_\pm \circ s_\pm^{-1})^* \frac{\partial}{\partial \alpha_i}[\omega]= (i_\pm \circ s_\pm^{-1})^* c_i, 
$$
since $\alpha_i^\pm=\alpha_i$ for $1\leq i \leq m-2$.
Moreover,
$$
c_{m-1}^+ =  \frac{\partial}{\partial \alpha^+_{m-1}}[\omega^+]= \left(\frac{\partial}{\partial \alpha_{m-1}}[\omega^+]\right)\left( \frac{\partial \alpha_{m-1}}{\partial \alpha_{m-1}^+} \right) = \frac{\partial}{\partial \alpha_{m-1}}[\omega^+] = (i_+ \circ s_+^{-1})^* c_{m-1},
$$
since $\alpha_{m-1}^+=\alpha_{m-1}+\alpha_m$. Similarly,
\begin{align*}
c_{m-1}^-  & =  \frac{\partial}{\partial \alpha^-_{m-1}}[\omega^-]= \left(\frac{\partial}{\partial \alpha_{m-1}}[\omega^-]\right)\left( \frac{\partial \alpha_{m-1}}{\partial \alpha_{m-1}^-} \right)= \frac{\partial}{\partial \alpha_{m-1}}[\omega^-] \cdot \,  \sgn(\alpha_{m-1}-\alpha_m)  \\  & =  \sgn(\alpha_{m-1}-\alpha_m) \, (i_- \circ s_-^{-1})^* c_{m-1},
\end{align*}
since $\alpha_{m-1}^-=\vert \alpha_{m-1}- \alpha_m\vert = \sgn( \alpha_{m-1}- \alpha_m) (\alpha_{m-1}- \alpha_m)$.
On the other hand, 
$$
c_{m-1}^+ =  \frac{\partial}{\partial \alpha^+_{m-1}}[\omega^+]= \left(\frac{\partial}{\partial \alpha_{m}}[\omega^+]\right)\left( \frac{\partial \alpha_{m}}{\partial \alpha_{m-1}^+} \right) = \frac{\partial}{\partial \alpha_{m}}[\omega^+] = (i_+ \circ s_+^{-1})^* c_{m}
$$
and 
\begin{align*}
c_{m-1}^-  & =  \frac{\partial}{\partial \alpha^-_{m-1}}[\omega^-]= \left(\frac{\partial}{\partial \alpha_{m}}[\omega^-]\right)\left( \frac{\partial \alpha_{m}}{\partial \alpha_{m-1}^-} \right)= - \frac{\partial}{\partial \alpha_{m}}[\omega^-] \cdot \,  \sgn(\alpha_{m-1}-\alpha_m)  \\  & = - \, \sgn(\alpha_{m-1}-\alpha_m) \, (i_- \circ s_-^{-1})^* c_{m}.
\end{align*}
\endproof
Finally, using Propositions~\ref{prop:4.1} and \ref{prop:5.1} we have our recursion formula.
\begin{Theorem}
\label{thm:5.1} 
Let  $\alpha^\pm$ and $c_i^\pm$  be as in Proposition~\ref{prop:5.1}. Then, for $k_1,\ldots, k_m \in \Z_{\geq 0}$ such that $k_1+\cdots + k_m= m-3$ and $k_m\geq 1$,
\begin{equation}\label{eq:5.1}
\begin{split}
&\dint\limits_{\M(\alpha)} c_1^{k_1} \cdots c_m^{k_m} = 
\dint\limits_{\M(\alpha^+)}
 (c^+_1)^{k_1} \cdots (c^+_{m-2})^{k_{m-2}} (c_{m-1}^+)^{k_{m-1}+k_m - 1} \quad + \\
  \\
& + \,  (-1)^{k_m - 1} \left( \sgn( \alpha_{m-1} -\alpha_m ) \right)^{k_{m-1}+k_m}
 \dint\limits_{\M(\alpha^-)}\!\! (c^-_1)^{k_1}\cdots
(c^-_{m-2})^{k_{m-2}} (c^-_{m-1})^{k_{m-1}+k_m - 1}. 
\end{split}
\end{equation}
\end{Theorem}  

\section{An explicit formula}\label{se:6}
Using this recursion formula we may obtain an explicit expression for the computation of intersection numbers. For that we first need to introduce  the following definition.
\begin{Definition}
Let $\alpha=(\alpha_1,\ldots, \alpha_m)$ be generic. A set $J\subset I := \{3,\ldots, m \}$ is called \emph{triangular} if 
$$
\sum_{i=3}^m (-1)^{\chi_{I\setminus J}(i)} \, \alpha_i > 0
$$
(where, for a set $S$, $\chi_S: S\to \{0,1\}$ is the characteristic function of $S$),
and  satisfies the following triangle inequalities:
\begin{align*}
&\alpha_1  \leq \alpha_2 + \sum_{i=3}^m (-1)^{\chi_{I\setminus J}(i)} \, \alpha_i, \\
&\alpha_2  \leq \alpha_1 + \sum_{i=3}^m (-1)^{\chi_{I\setminus J}(i)} \, \alpha_i, \\
& \sum_{i=3}^m (-1)^{\chi_{I\setminus J}(i)} \, \alpha_i  \leq \alpha_1 + \alpha_2.
\end{align*}
Moreover, define
$\mathcal{T}:=\mathcal{T}(\alpha) :=\{ J\subset I\mid \, J \, \text{is triangular}\}$, the family of all triangular sets in $I$.
\end{Definition} 
Finally we come to Theorem~\ref{thm:1.2}.
\begin{Theorem}\label{thm:6.1}
Let $\alpha=(\alpha_1,\ldots, \alpha_m)$ be generic. Suppose $k_{m-l},\ldots, k_m\in \Z_+$,  $k_1=\cdots = k_{m-l-1}=0$ and $k_{m-l} + \cdots + k_m = m-3$. Let  $c_i:=c_1(V_i(\alpha))$ be the first Chern classes of the circle bundles $V_i(\alpha)\to \M(\alpha)$. Then
\begin{equation}
\label{eq:6.1}
\int_{\M(\alpha)} c_{m-l}^{k_{m-l}} \cdots c_m^{k_m} = \sum_{J \in \mathcal{T}(\alpha)} (-1)^{\left( \sum_{i\in I\setminus J} k_i \right) + m - \vert J  \vert}.
\end{equation}
\proof
We will prove this formula by induction on $m$ starting with $m=4$.

For $m=4$, the recursion formula of Theorem~\ref{thm:5.1} gives
$$
\int_{\M(\alpha_1,\ldots,\alpha_4)} c_4 = \int_{\M(\alpha_1,\alpha_2,\alpha_3+\alpha_4)} 1 +\, \sgn(\alpha_3-\alpha_4)\int_{\M(\alpha_1,\alpha_2,\vert \alpha_3 - \alpha_4 \vert)} 1.
$$ 
Hence,
$$
\int_{\M(\alpha_1,\ldots,\alpha_4)} c_4 = 
\left\{ \begin{array}{cll} 
2 & & \text{if} \quad \mathcal{T}=\left\{\{3\},\{3,4\}\right\} \\
0 & & \text{if} \quad \mathcal{T}=\left\{\{4\},\{3,4\}\right\} \\ 
1 & & \text{if} \quad \mathcal{T}=\left\{\{3\}\right\} \, \text{or} \, \left\{\{3,4\}\right\} \\ 
-1 & & \text{if} \quad \mathcal{T}=\left\{\{4\}\right\},
\end{array} \right.
$$
and it is easy to verify that this agrees in all cases with  the right-hand side of \eqref{eq:6.1}.

We will now assume that  \eqref{eq:6.1} holds for some $m$ and show that it is still valid for $m+1$. Using the recursion formula \eqref{eq:5.1} once more we obtain
\begin{align*}
 & \dint \limits_{\M(\alpha_1,\ldots, \alpha_m,\alpha_{m+1})} \hspace{-1cm} c_{m+1-l}^{k_{m+1-l}} \cdots c_{m+1}^{k_{m+1}} \, = \hspace{-1.2cm} \int \limits_{\M(\alpha_1,\ldots,\alpha_{m-1},\alpha_m+\alpha_{m+1})}\hspace{-1.5cm} (c_{m+1-l}^+)^{k_{m+1-l}} \cdots (c_{m-1}^+)^{k_{m-1}} (c_{m}^+)^{k_{m}+k_{m+1}-1} \, + \\  \\ & + (-1)^{k_{m+1}-1}\,\sgn(\alpha_m - \alpha_{m+1})^{k_m + k_{m+1}}\hspace{-1.8cm} \int \limits_{\M(\alpha_1,\ldots,\alpha_{m-1},\vert \alpha_m - \alpha_{m+1}\vert)} \hspace{-1.5cm} (c_{m+1-l}^-)^{k_{m+1-l}} \cdots(c_{m-1}^-)^{k_{m-1}} (c_{m}^-)^{k_{m}+k_{m+1}-1}. \\ 
\end{align*}
Writing 
\begin{align*}
\mathcal{T}_{m+1}: = \mathcal{T}(\alpha_1,\ldots,\alpha_m, \alpha_{m+1}), & \quad  \mathcal{T}_m^+ := \mathcal{T}(\alpha_1,\ldots,\alpha_{m-1}, \alpha_m + \alpha_{m+1}) \quad \text{and} \\  \mathcal{T}_m^- = \mathcal{T}&(\alpha_1,\ldots,\alpha_{m-1}, \vert \alpha_m - \alpha_{m+1}\vert ),
\end{align*} 
then, if $\alpha_{m}-\alpha_{m+1} > 0$,
\begin{align*}
\mathcal{T}_{m+1} = & \left\{\widetilde{J} \in \mathcal{T}_m^+ \mid m \notin \widetilde{J} \right\} \,  \bigcup \, \left\{\widetilde{J} \cup \{m+1\}\mid \widetilde{J} \in \mathcal{T}_m^+ \, \text{and} \, m \in \widetilde{J} \right\} \, \bigcup \\ & \bigcup \,  \left\{\widetilde{J} \in \mathcal{T}_m^- \mid m \in \widetilde{J} \right\} \, \bigcup \,   \left\{\widetilde{J}\cup \{m+1\} \mid \widetilde{J} \in \mathcal{T}_m^- \, \text{and} \, m \notin \widetilde{J} \right\},
\end{align*}
while, if $\alpha_{m}-\alpha_{m+1} < 0$,
\begin{align*}
& \mathcal{T}_{m+1} =  \left\{\widetilde{J} \in \mathcal{T}_m^+ \mid m \notin \widetilde{J} \right\} \,  \bigcup \, \left\{\widetilde{J} \cup \{m+1\}\mid \widetilde{J} \in \mathcal{T}_m^+ \, \text{and} \, m \in \widetilde{J} \right\} \, \bigcup \\  \bigcup \, &  \left\{\left(\widetilde{J}\setminus\{m\}\right)\cup\{m+1\}\mid \widetilde{J} \in \mathcal{T}_m^- \, \text{and} \, m \in \widetilde{J} \right\} \, \bigcup \,   \left\{\widetilde{J}\cup \{m\} \mid \widetilde{J} \in \mathcal{T}_m^- \, \text{and} \, m \notin \widetilde{J} \right\}.
\end{align*}
Assuming  \eqref{eq:6.1} holds and writing  $\tilde{k}_j=k_j$ for $j\neq m$ and $\tilde{k}_m=k_m + k_{m+1} -1$, we have
\begin{align*}
 \int \limits_{\M(\alpha_1,\ldots,\alpha_{m-1},\alpha_m+\alpha_{m+1})}\hspace{-1.5cm}  &(c_{m+1-l}^+)^{k_{m+1-l}}\cdots   (c_{m-1}^+)^{k_{m-1}} (c_{m}^+)^{k_{m}+k_{m+1}-1}  \\ & = \sum_{\widetilde{J}\in \mathcal{T}_m^+} (-1)^{\left(\sum_{i\in \{3,\ldots,m\}\setminus \widetilde{J}} \tilde{k}_i\right)+m-\vert \widetilde{J} \vert} \\ & = \sum_{\widetilde{J}\in \mathcal{T}_m^+\,\text{s.t.} m\notin \widetilde{J}} (-1)^{\left(\sum_{i\in \{3,\ldots,m-1\}\setminus \widetilde{J}} \, k_i\right) + \, \tilde{k}_m + \, m -\vert \widetilde{J} \, \vert} \,\, + \\ & +  \sum_{\widetilde{J}\in \mathcal{T}_m^+\,\text{s.t.} m \in \widetilde{J}} (-1)^{\left(\sum_{i\in \{3,\ldots,m+1\}\setminus (\widetilde{J} \, \cup \, \{m+1\})} k_i\right)+ \, m - \vert \widetilde{J} \, \vert} \\ & = \sum_{\widetilde{J}\in \mathcal{T}_m^+\,\text{s.t.} m \notin \widetilde{J}} (-1)^{\left(\sum_{i \in \{3,\ldots,m+1\}\setminus \widetilde{J}}\, k_i\right)+  m-1 -\vert \widetilde{J} \vert} \,\, + \\ & +  \sum_{\widetilde{J}\in \mathcal{T}_m^+\,\text{s.t.} m\in \widetilde{J}} (-1)^{\left(\sum_{i\in \{3,\ldots,m+1\}\setminus (\widetilde{J} \, \cup \, \{m+1\})} k_i\right)+ \, m + 1 -\vert \widetilde{J}\, \cup \, \{m+1\} \vert}.  
\end{align*}
On the other hand, if $\alpha_m - \alpha_{m+1}>0$,
\begin{align*}
& (-1)^{k_{m+1}-1}\,\,\sgn(\alpha_m - \alpha_{m+1})^{k_m + k_{m+1}}\hspace{-1cm} \int \limits_{\M(\alpha_1,\ldots,\alpha_{m-1},\vert \alpha_m - \alpha_{m+1}\vert)} \hspace{-1.5cm}   (c_{m+1-l}^-)^{k_{m+1-l}} \cdots (c_{m-1}^-)^{k_{m-1}} (c_{m}^-)^{k_{m}+k_{m+1}-1} \\ & = (-1)^{k_{m+1}-1}\,\,\sum_{\widetilde{J}\in \mathcal{T}_m^-} (-1)^{\left(\sum_{i\in \{3,\ldots,m\}\setminus \widetilde{J}} \tilde{k}_i\right) + m - \vert \widetilde{J} \vert} \\ & =
 \sum_{\widetilde{J}\in \mathcal{T}_m^-\,\text{s.t.} m \in \widetilde{J}} (-1)^{\left(\sum_{i\in \{3,\ldots,m\}\setminus \widetilde{J}} \, k_i\right) + \, (k_{m+1}-1) + \, m -\vert \widetilde{J} \, \vert} \\ & +  \sum_{\widetilde{J}\in \mathcal{T}_m^-\,\text{s.t.} m \notin \widetilde{J}} (-1)^{\left(\sum_{i\in \{3,\ldots,m-1\}\setminus \widetilde{J}} k_i \right)+ \overbrace{k_m +k_{m+1} -1}^{\tilde{k}_m} + \, (k_{m+1}-1) + \, m - \vert \widetilde{J} \, \vert} \\ & = \sum_{\widetilde{J}\in \mathcal{T}_m^-\,\text{s.t.} m \in \widetilde{J}} (-1)^{\left(\sum_{i \in \{3,\ldots,m+1\}\setminus \widetilde{J}}\, k_i\right)+  m-1 -\vert \widetilde{J} \vert} \\ & +  \sum_{\widetilde{J}\in \mathcal{T}_m^-\,\text{s.t.} m\notin \widetilde{J}} (-1)^{\left(\sum_{i\in \{3,\ldots,m+1\}\setminus (\widetilde{J} \, \cup \, \{m+1\})} k_i\right)+ \, m + 1 -\vert \widetilde{J}\, \cup \, \{m+1\} \vert},
\end{align*}
while,  if $\alpha_m - \alpha_{m+1}<0$,
\begin{align*}
& (-1)^{k_{m+1}-1}\,\,\sgn(\alpha_m - \alpha_{m+1})^{k_m + k_{m+1}}\hspace{-1cm} \int \limits_{\M(\alpha_1,\ldots,\alpha_{m-1},\vert \alpha_m - \alpha_{m+1}\vert)} \hspace{-1.5cm}   (c_{m+1-l}^-)^{k_{m+1-l}} \cdots (c_{m-1}^-)^{k_{m-1}}(c_{m}^-)^{k_{m}+k_{m+1}-1} \\ & = (-1)^{k_{m+1}-1+ k_m + k_{m+1}}\, \sum_{\widetilde{J}\in \mathcal{T}_m^-} (-1)^{\left(\sum_{i\in \{3,\ldots,m\}\setminus \widetilde{J}} \tilde{k}_i\right) + m - \vert \widetilde{J} \vert} \\ & = 
 \sum_{\widetilde{J}\in \mathcal{T}_m^-\,\text{s.t.} m \in \widetilde{J}} (-1)^{\left(\sum_{i\in \{3,\ldots,m \}\setminus \widetilde{J}} \, k_i\right) + k_m -1 + \, m -\vert \widetilde{J} \, \vert} \\ & +  \sum_{\widetilde{J}\in \mathcal{T}_m^-\,\text{s.t.} m \notin \widetilde{J}} (-1)^{\left(\sum_{i\in \{3,\ldots,m-1\}\setminus \widetilde{J}} k_i \right)+ \overbrace{k_m + k_{m+1} - 1}^{\tilde{k}_m} + k_m - 1+ \, m - \vert \widetilde{J} \, \vert} \\ & = \sum_{\widetilde{J}\in \mathcal{T}_m^-\,\text{s.t.} m \in \widetilde{J}} (-1)^{\left(\sum_{i \in \{3,\ldots,m+1\}\setminus ((\widetilde{J}\setminus\{m\})\cup \{m+1\})}\, k_i\right)+  m-1 -\vert ( \widetilde{J}\setminus\{m\})\cup\{m+1\} \vert} \\ & +  \sum_{\widetilde{J}\in \mathcal{T}_m^-\,\text{s.t.} m\notin \widetilde{J}} (-1)^{\left(\sum_{i\in \{3,\ldots,m+1\}\setminus (\widetilde{J} \, \cup \, \{m\})} k_i\right)+ \, m + 1 -\vert \widetilde{J}\, \cup \, \{m\} \vert} 
\end{align*}
and the result follows.

Note that in the above proof we have to assume that each time we use the recursion formula we have $\alpha_m\neq \alpha_{m+1}$. Nevertheless, the result is still valid even if this is not the case, as long as $\alpha$ is generic.  Indeed, for a generic  $\alpha$ such that $\alpha_{m} = \alpha_{m+1}$ we may use Remark~\ref{rk:2.1} and take a small value of $\varepsilon > 0$ for which
$\M(\alpha)$ is diffeomorphic to $\M(\alpha_\varepsilon)$ with $\alpha:=(\alpha_1,\ldots,\alpha_{m}, \alpha_{m+1})$ and $\alpha_\varepsilon:=(\alpha_1,\ldots,\alpha_{m}, \alpha_{m} + \varepsilon)$. Note that, for $\varepsilon$ small enough, we have $\mathcal{T}^\pm_m(\alpha)=\mathcal{T}^\pm_m(\alpha_\varepsilon)$ (since $\alpha$ generic implies that $\alpha^+$ and $\alpha^-$ are also generic) and so the induction step still holds. 
\endproof
\end{Theorem}

\section{Examples}\label{se:7}
\begin{Example}
\label{ex:7.1}
Let us consider $\M(4,3,4,3,4)$.  It is a toric manifold of dimension $4$  obtained by symplectic blowing up three points in $S^2\times S^2$. This can  be seen for instance in the moment polytope depicted in Figure~\ref{fig:7.1} obtained from the Hamiltonian $2$-torus action given by the  bending flows along the second and third diagonals. This polytope is the intersection of the rectangle $[1,7]\times [1,7]$ with the non-compact rectangular region
$$
\{(x,y)\in (\R_{\geq 0})^2\mid\,  x + y \geq 4, \, y \geq x - 4 \, \text{and} \,\, y\leq x + 4 \}
$$
(cf. \cite{HK2} for details on how to obtain these moment polytopes).
\begin{figure}[h!]
\begin{center}
\scalebox{.8}{\includegraphics*[0mm,210mm][80mm,250mm]{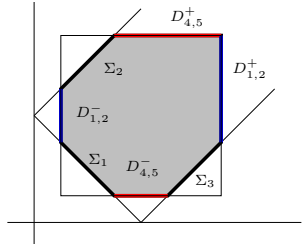}}
\end{center}
\caption{Moment map image for the bending action along the second and third diagonals of $\M(4,3,4,3,4)$.}\label{fig:7.1}
\end{figure}

 We explicitly compute some intersection numbers associated to this manifold. By Remark~\ref{rk:3.1} we  note that
\begin{equation}\nonumber
\dint\limits_{\M(4,3,4,3,4)} \!\!\!\!\!\!\!\!\!\! c_1^2 \,\, = \dint\limits_{\M(4,3,4,3,4)} \!\!\!\!\!\!\!\!\!\! c_3^2 \,\, = \dint\limits_{\M(4,3,4,3,4)} \!\!\!\!\!\!\!\!\!\! c_5^2.
\end{equation}
Then, by Theorem~\ref{thm:5.1}, we have
\begin{align*}
&\dint\limits_{\M(4,3,4,3,4)} \!\!\!\!\!\!\!\!\!\! c_5^2 \, \, = \dint\limits_{\M(4,3,4,7)} \!\!\!\!\!\!\!\!\!\! c^+_4 \,\, \begin{array}{c} \\ - \\ \fbox{$\scriptscriptstyle{\alpha_4 < \alpha_5}$} \end{array}
\dint\limits_{\M(4,3,4,1)}
\!\!\!\!\!\!\!\!\!\! c^-_4 =  \\
& = \left( \dint\limits_{\M(4,3,11)} 1 \, \begin{array}{c} \\ - \\ \fbox{$\scriptscriptstyle{\alpha^+_3<\alpha^+_4}$} \end{array} \dint\limits_{\M(4,3,3)} 1 \right)  - \left( \dint\limits_{\M(4,3,5)} 1 \begin{array}{c} \\ + \\ \fbox{$\scriptscriptstyle{\alpha^-_3 > \alpha^-_4}$} \end{array}  \dint\limits_{\M(4,3,3)} 1  \right)  \\
& = 0 - 1 - 1 - 1 = -3 .
\end{align*}
Note that $\M(4,3,11)=\varnothing$ since the triple $(4,3,11)$ does not satisfy the necessary triangle inequalities ($11> 4 + 3$). 
This result can be confirmed by taking the intersections of the corresponding manifolds $D_{i,j}$. For instance, since $\alpha_4 < \alpha_5$, the Poincar\'e dual of $c_5$ is $[D_{4,5}^+] - [D_{4,5}^-]$ (cf. Proposition~\ref{prop:4.1}), and so
\begin{align*}
& \int \limits_{\M(4,3,4,3,4)}  c_5^2 = \int_{\text{PD}(c_5)} c_5\vert_{\text{PD}(c_5)}= \text{PD}(c_5) \bullet  \text{PD}(c_5) = \\ & =   ([D_{4,5}^+] - [D_{4,5}^-]) \bullet ([D_{4,5}^+] - [D_{4,5}^-])   = [D_{4,5}^+] \bullet  [D_{4,5}^+] +  [D_{4,5}^-] \bullet  [D_{4,5}^-] = -1 - 2 =-3.  
\end{align*}
The image of the manifolds  $D_{4,5}^\pm$ by the moment map can be seen in Figure~\ref{fig:7.1}. Note that their intersection numbers can also be read from the moment polytope. They are negative and their absolute value equals the absolute value of the determinant of the $2\times 2$-matrix formed by the inward unit normal vectors to the adjacent edges. In particular, for $D^+_{4,5}$ and $
D^-_{4,5}$ we have 
$$
D^+_{4,5}\bullet D^+_{4,5} = - \left\vert \det \left( \begin{array}{rr} -1 & 1 \\ 0 & - 1 \end{array}\right)\right\vert = -1 \,\, \text{and} \,\, D^-_{4,5}\bullet D^-_{4,5} = - \left\vert \det \left( \begin{array}{rr} - 1 & 1 \\ 1 & 1 \end{array}\right) \right\vert = -2.
$$ 

Similarly, by Remark~\ref{rk:3.1},
\begin{equation}\nonumber
\dint\limits_{\M(4,3,4,3,4)} \!\!\!\!\!\!\!\!\!\! c_1c_2 =  \dint\limits_{\M(4,3,4,3,4)} \!\!\!\!\!\!\!\!\!\! c_1c_4 =
\dint\limits_{\M(4,3,4,3,4)} \!\!\!\!\!\!\!\!\!\! c_2c_3 =  \dint\limits_{\M(4,3,4,3,4)} \!\!\!\!\!\!\!\!\!\! c_2c_5 =
\dint\limits_{\M(4,3,4,3,4)} \!\!\!\!\!\!\!\!\!\! c_3c_4 =  \dint\limits_{\M(4,3,4,3,4)} \!\!\!\!\!\!\!\!\!\! c_4c_5,
\end{equation}
and then, by Theorem~\ref{thm:5.1},
\begin{align*}
& \dint\limits_{\M(4,3,4,3,4)} \!\!\!\!\!\!\!\!\!\! c_4c_5 = \dint\limits_{\M(4,3,4,7)} \!\!\!\!\!\!\!\!\!\! c_4^+ \, \begin{array}{c} \\ + \\ \fbox{$\scriptscriptstyle{\alpha_4 < \alpha_5}$} \end{array} \, \dint\limits_{\M(4,3,4,1)}  \!\!\!\!\!\!\!\!\!\! c_4^- \\ 
& = \left(\dint\limits_{\M(4,3,11)} 1  \begin{array}{c} \\  - \\ \fbox{$\scriptscriptstyle{\alpha^+_3 < \alpha^+_4}$} \end{array} \, \dint\limits_{\M(4,3,3)} 1 \right)
+ \left( \dint\limits_{\M(4,3,5)} 1 \,  \begin{array}{c} \\ +\\ \fbox{$\scriptscriptstyle{\alpha^-_3 > \alpha^-_4}$} \end{array}  \, \dint\limits_{\M(4,3,3)} 1 \, \right)  \\ & = 0 -1 +1 +1 =1.\\
\end{align*}
Moreover, by Remarks~\ref{rk:2.1} and \ref{rk:3.1},
\begin{align*}
&\dint\limits_{\M(4,3,4,3,4)} \!\!\!\!\!\!\!\!\!\! c_2c_4 =  \dint\limits_{\M(4,4,4,3,3+\varepsilon)} \!\!\!\!\!\!\!\!\!\! c_4c_5 = \dint\limits_{\M(4,4,4,6)} \!\!\!\!\!\!\!\!\!\! c_4^+ \, \begin{array}{c} \\ + \\ \fbox{$\scriptscriptstyle{\alpha_4 < \alpha_5}$} \end{array} \, \dint\limits_{\M(4,4,4,\varepsilon)}  \!\!\!\!\!\!\!\!\!\! c_4^- = \\ 
& = \left(\dint\limits_{\M(4,4,10)} 1  \begin{array}{c} \\  - \\ \fbox{$\scriptscriptstyle{\alpha^+_3 < \alpha^+_4}$} \end{array} \, \dint\limits_{\M(4,4,2)} 1 \right)
+ \left( \dint\limits_{\M(4,4,4+\varepsilon)} 1 \,  \begin{array}{c} \\ +\\ \fbox{$\scriptscriptstyle{\alpha^-_3 > \alpha^-_4}$} \end{array}  \, \dint\limits_{\M(4,4,4-\varepsilon)} 1 \, \right)  \\ & = 0 -1 +1 +1 = 1.
\end{align*}
Again, these results can be obtained by taking the Poincar\'e duals. For example,
\begin{align*}
& \!\!\!\!\!\!\!\!\! \dint\limits_{\M(4,3,4,3,4)} \!\!\!\!\!\!\!\!\!\! c_2c_4 = \text{PD}(c_2) \bullet  \text{PD}(c_4) = (D_{1,2}^+ + D_{1,2}^-)\bullet (D_{4,5}^+ + D_{4,5}^-) =  D_{1,2}^+ \bullet D_{4,5}^+ = 1.
\end{align*}
If instead of the recursion formula \eqref{eq:5.1} we use Theorem~\ref{thm:6.1} we get for instance 
\begin{align*}
 \dint\limits_{\M(4,3,4,3,4)} \!\!\!\!\!\!\!\!\!\! c_5^2 \, \, & = \sum_{J \in \mathcal{T}(4,3,4,3,4)} (-1)^{\left( \sum_{i \in I\setminus J} k_i \right) + m - \vert J  \vert} \\ \\ & = (-1)^{k_5+3} + (-1)^{k_4+3} + (-1)^{k_3+3}  = -3,
\end{align*}
since $\mathcal{T}(4,3,4,3,4)=\{ \{3,4\},\{3,5\},\{4,5\}\}$, $k_5=2$, $k_4=k_3=0$ and $m-\vert J \vert=3$. Moreover, 
\begin{align*}
\dint\limits_{\M(4,3,4,3,4)} \!\!\!\!\!\!\!\!\!\! c_2^2 \, \, & =  \dint\limits_{\M(4,4,4,3,3)} c_5^2 = \dint\limits_{\M(4,4,4,3,3+\varepsilon)} \!\!\!\!\!\!\!\!\!\! c_5^2 = \sum_{\tiny{ \begin{array}{c} J \in \mathcal{T}(4,4,4,3,3) \\ = \mathcal{T}(4,4,4,3,3+\varepsilon) \end{array} }} (-1)^{\left( \sum_{i \in I\setminus J} k_i \right) + m - \vert J  \vert} \\ \\ & = (-1)^{k_5+3} + (-1)^{k_4+3} + (-1)^{k_3+3}  = -3,
\end{align*}
since $\mathcal{T}(4,4,4,3,3)=\{ \{3,4\},\{3,5\},\{4,5\}\}$, $k_5=2$ and $k_4=k_3=0$. 
               
Using similar arguments we can easily obtain the complete list of intersection numbers for $\M(4,3,4,3,4)$ (cf. Table~\ref{tab:7.1}).

\begin{table}[h!]
\begin{center}
\begin{tabular}{||c|c||}
\hline
$\int\limits_{M(\alpha)}c_i^2=-3$   & $i=1,\ldots, 5$ \\
\hline
$\int\limits_{M(\alpha)}c_ic_j=1$   & ($i\neq j$), $i,j=1, \ldots, 5$\\ 
\hline
\end{tabular}
\end{center}
\caption{Intersection numbers for $\M(4,3,4,3,4)$.}
\label{tab:7.1}
\end{table}

The exceptional divisors $\Sigma_1$, $\Sigma_2$ and $\Sigma_3$ of the blow ups of $S^2\times S^2$ are:
\begin{align*}
\Sigma_1 =  \frac{1}{2}\left( D_{1,2}^+ - D_{1,2}^-  + D_{4,5}^+ - D_{4,5}^- \right)&,\quad 
\Sigma_2  =  \frac{1}{2}\left( D_{1,2}^+ - D_{1,2}^- + D_{2,3}^+ - D_{2,3}^- \right) \quad \text{and} \quad \\
\Sigma_3   = \frac{1}{2}& \left(  D_{2,3}^+ - D_{2,3}^- + D_{4,5}^+ - D_{4,5}^-\right),
\end{align*}
and so 
\begin{align*}
\text{PD}(\Sigma_1)  =  \frac{1}{2} (c_1 + c_5) ,\quad 
\text{PD}(\Sigma_2)  =  \frac{1}{2} (c_1+c_3) \quad \text{and} \quad 
\text{PD}(\Sigma_3)  =  \frac{1}{2} (c_3+c_5). 
\end{align*}
Note incidentally that $\int_{\M(\alpha)} c_i^2$ has the same value for every $i$. This agrees with the fact proved in \cite{HK2} that in any polygon space we have $c_i^2=p$ for all $i$, where $p$ is the Pontrjagin class of the principal $SO(3)$-bundle $A(\alpha) \to \M(\alpha)$, with 
$$
A(\alpha):= \{ (\vec{v}_1,\ldots,\vec{v}_m)\in S^2_{\alpha_i}\mid \sum \vec{v}_i =0\}.
$$
\end{Example}

\begin{Example}
Let $\alpha=(1,\ldots,1,m-2)$. We know from \cite{H} that $\M(\alpha)$ is diffeomorphic to $\C P^{m-3}$ so we will check our results.

For nonnegative integers $k_1,\ldots,k_m$ with $k_1 + \cdots + k_m = m-3$ we obtain
\begin{equation}\nonumber
\dint\limits_{\M(1,\ldots,1,m-2)} \!\!\!\!\!\!\!\!\!\! c_1^{k_1} \cdots c_{m}^{k_m} \,\, = (-1)^{k_m}.
\end{equation}
Indeed, if $k_m\geq 1$, we have by Remark~\ref{rk:3.1} that
\begin{align*}
\dint\limits_{\M(1,\ldots,1,m-2)} \!\!\!\!\!\!\!\!\!\! c_1^{k_1} \cdots c_{m}^{k_m} \,\, & =  \dint\limits_{\M(1,\ldots,1,m-2)} c_{m-l}^{\tilde{k}_{m-l}} \cdots  c_{m-1}^{\tilde{k}_{m-1}}  c_{m}^{k_m}  
\end{align*}
with $\tilde{k}_{m-l}, \ldots , \tilde{k}_{m-1} \geq 1$, $m-l > 3$ and  $\tilde{k}_{m-l}+ \cdots + \tilde{k}_{m-1}+k_m=m-3$. Then, since $\mathcal{T} (1,\ldots,1,m-2) =\left\{\{m\}\right\}$, Theorem~\ref{thm:6.1} yields
\begin{equation}\nonumber
\dint\limits_{\M(1,\ldots,1,m-2)} \!\!\!\!\!\!\!\!\!\! c_1^{k_1} \cdots c_{m}^{k_m} \,\, = (-1)^{\tilde{k}_{m-l}+\cdots+\tilde{k}_{m-1} + m -1} = (-1)^{m - 3- k_m + m -1}=(-1)^{k_m}.
\end{equation}
On the other hand, if $k_m=0$, again using Remark~\ref{rk:3.1}, we get
\begin{align*}
\dint\limits_{\M(1,\ldots,1,m-2)} \!\!\!\!\!\!\!\!\!\! c_1^{k_1} \cdots c_{m}^{k_m} \,\, & =  \!\!\!\!\!\!\!\!\!\! \dint\limits_{\M(1,\ldots,1,m-2)}  \!\!\!\!\!\!\!\!\!\! c_{m-l}^{\tilde{k}_{m-l}} \cdots  c_{m-1}^{\tilde{k}_{m-1}}  =  \!\!\!\!\!\!\!\!\!\! \!\!\!\!\!\!\!\!\!\! \dint\limits_{\M(1,\ldots,1,\underbrace{m-2}_{m-l},1,\ldots,1)} \!\!\!\!\!\!\!\!\!\! \!\!\!\!\!\!\!\!\!\! c_{m-l+1}^{\tilde{k}_{m-l+1}} \cdots  c_{m-1}^{\tilde{k}_{m-1}} c_{m}^{\tilde{k}_{m-l}},  
\end{align*}
with  $\tilde{k}_{m-l}, \ldots, \tilde{k}_{m-1} \geq 1$, $m-l \geq 3$ and  $\tilde{k}_{m-l} + \cdots + \tilde{k}_{m-1} = m-3$. Then, since $\mathcal{T}(1,\ldots,1,m-2,1,\ldots,1)=\left\{\{m-l\}\right\}$, Theorem~\ref{thm:6.1} gives
\begin{equation}\nonumber
\dint\limits_{\M(1,\ldots,1,m-2)} \!\!\!\!\!\!\!\!\!\! c_1^{k_1} \cdots c_{m}^{k_m} \,\, = (-1)^{\tilde{k}_{m-l+1}+\cdots+\tilde{k}_{m-1} + \tilde{k}_{m-l} + m -1} = (-1)^{m - 3 + m -1} = 1= (-1)^{k_m}.
\end{equation}  
We conclude that $c_1=\cdots=c_{m-1}=h$ where $h\in H^2(\C P^{m-3}; \Z)$ is the canonical generator, while $c_m=-h$.
\end{Example}

\begin{Example}
Let $\alpha=(\varepsilon,\ldots,\varepsilon,1,1,1)$ with $(m-3)\varepsilon <1$. We know from \cite{HK2} that $\M(\alpha)$ is symplectomorphic to $\prod_{i=1}^{m-3} S^2_{\alpha_i}$.

For nonnegative integers $k_1,\ldots,k_m$ with $k_1 + \cdots + k_m = m-3$ we have
\begin{equation}\label{eq:7.1}
\dint\limits_{\M(\varepsilon, \ldots, \varepsilon, 1,1,1)} \!\!\!\!\!\!\!\!\!\! c_1^{k_1} \cdots c_{m}^{k_m} \,\, = 0
\end{equation}
unless $k_{m-2}=k_{m-1}=k_m=0$ and $k_1=\cdots =k_{m-3}=1$, in which case we have
\begin{equation}\nonumber
\dint\limits_{\M(\varepsilon, \ldots, \varepsilon, 1,1,1)} \!\!\!\!\!\!\!\!\!\! c_1 \cdots c_{m-3} \,\, = 2^{m-3}.
\end{equation}
Indeed, if  $k_{m-2},k_{m-1},k_m\neq 0$ then  by Remark~\ref{rk:3.1}, 
\begin{align*}
\dint\limits_{\M(\varepsilon, \ldots, \varepsilon, 1,1,1)} \!\!\!\!\!\!\!\!\!\! c_1^{k_1} \cdots c_{m}^{k_m} \,\, = \dint\limits_{\M(\varepsilon, \ldots, \varepsilon, 1,1,1)} \!\!\!\!\!\!\!\!\!\! c_{m-l}^{\tilde{k}_{m-l}} \cdots  c_{m-3}^{\tilde{k}_{m-3}}  c_{m-2}^{k_{m-2}} c_{m-1}^{k_{m-1}} c_{m}^{k_m} \,\,
\end{align*}
with $\tilde{k}_{m-l}, \ldots, \tilde{k}_{m-3} \geq 1$ and $m-l > 3$. Then, \eqref{eq:7.1} holds since $\mathcal{T}(\varepsilon, \ldots, \varepsilon,1,1,1)=\varnothing$.

If only one of $k_{m-2}$, $k_{m-1}$ or $k_m$ is equal to zero then, again by Remark~\ref{rk:3.1},
\begin{align*}
\dint\limits_{\M(\varepsilon, \ldots, \varepsilon, 1,1,1)} \!\!\!\!\!\!\!\!\!\! c_1^{k_1} \cdots c_{m}^{k_m} \,\, = \dint\limits_{\M(\varepsilon, \ldots, \varepsilon, 1,1,1)} \!\!\!\!\!\!\!\!\!\! c_{m-l}^{\tilde{k}_{m-l}} \cdots c_{m-1}^{\tilde{k}_{m-1}} \,\,
\end{align*}
with $\tilde{k}_{m-l}, \ldots, \tilde{k}_{m-1} \geq 1$ and $m-l \geq 3$ and so
\begin{align*}
\dint\limits_{\M(\varepsilon, \ldots, \varepsilon, 1,1,1)} \!\!\!\!\!\!\!\!\!\! c_1^{k_1} \cdots c_{m}^{k_m} \,\, =    \dint\limits_{\M(\varepsilon, \ldots, \varepsilon, \underbrace{1}_{m-l}, \varepsilon, \ldots, \varepsilon, 1,1,\varepsilon)} \!\!\!\!\!\!\!\!\!\! c_{m-l+1}^{\tilde{k}_{m-l+1}} \cdots c_{m-1}^{\tilde{k}_{m-1}}c_m^{\tilde{k}_{m-l}} \,\,=0,
\end{align*}
since, for $m-l>3$,  $\mathcal{T}(\varepsilon, \ldots, \varepsilon, \underbrace{1}_{m-l}, \varepsilon, \ldots, \varepsilon, 1,1,\varepsilon)=\varnothing$. 

If two of $k_{m-2}$, $k_{m-1}$ or $k_m$ are zero then  by Remark~\ref{rk:3.1},
\begin{align*}
\dint\limits_{\M(\varepsilon, \ldots, \varepsilon, 1,1,1)} \!\!\!\!\!\!\!\!\!\! c_1^{k_1} \cdots c_{m}^{k_m} \,\, = \dint\limits_{\M(\varepsilon, \ldots, \varepsilon, 1,1,1)} \!\!\!\!\!\!\!\!\!\! c_{m-l}^{\tilde{k}_{m-l}} \cdots c_{m-2}^{\tilde{k}_{m-2}} \,\,
\end{align*}
with $\tilde{k}_{m-l}, \ldots, \tilde{k}_{m-2} \geq 1$ and $m-l \geq 2$ and so
\begin{align*}
\dint\limits_{\M(\varepsilon, \ldots, \varepsilon, 1,1,1)} \!\!\!\!\!\!\!\!\!\! c_1^{k_1} \cdots c_{m}^{k_m} \,\, 
=    \dint\limits_{\M(\varepsilon, \ldots, \varepsilon, \underbrace{1}_{m-l},1, \varepsilon, \ldots, \varepsilon, 1,\varepsilon,\varepsilon)} \!\!\!\!\!\!\!\!\!\! c_{m-l+2}^{\tilde{k}_{m-l+2}} \cdots c_{m-2}^{\tilde{k}_{m-2}} c_{m-1}^{\tilde{k}_{m-l}} c_m^{\tilde{k}_{m-l+1}} \,\,=0
\end{align*}
since both  $\mathcal{T} (\varepsilon, 1,1,\varepsilon \ldots, \varepsilon, 1,\varepsilon,\varepsilon )$ and  $\mathcal{T} (\varepsilon,\varepsilon, \ldots, \varepsilon, \underbrace{1}_{m-l},1, \varepsilon, \ldots, \varepsilon, 1,\varepsilon,\varepsilon)$ (respectively the cases  $m-l=2$ and $m-l>2$) are empty.

If $k_{m-2}=k_{m-1}=k_m=0$ then 
\begin{align*}
\dint\limits_{\M (\varepsilon, \ldots, \varepsilon, 1,1,1)} \!\!\!\!\!\!\!\!\!\! c_1^{k_1} \cdots c_{m}^{k_m} \,\, = \dint\limits_{\M (\varepsilon, \ldots, \varepsilon, 1,1,1)} \!\!\!\!\!\!\!\!\!\! c_{m-l}^{\tilde{k}_{m-l}} \cdots c_{m-3}^{\tilde{k}_{m-3}} \,\,
\end{align*}
with $\tilde{k}_{m-l}, \ldots, \tilde{k}_{m-3} \geq 1$ and $m-l \geq 1$. If $m-l\geq 2$ then we do as before and again obtain $\mathcal{T}=\varnothing$ implying that \eqref{eq:7.1} holds. If, however, $m-l=1$ then necessarily $k_1=\cdots = k_{m-3}=1$ and then
\begin{align*}
\dint\limits_{\M(\varepsilon, \ldots, \varepsilon, 1,1,1)} \!\!\!\!\!\!\!\!\!\! c_1^{k_1} \cdots c_{m}^{k_m} \,\, = \dint\limits_{\M(\varepsilon, \ldots, \varepsilon, 1,1,1)} \!\!\!\!\!\!\!\!\!\! c_1 \cdots c_{m-3} = \dint\limits_{\M(1,1,1, \varepsilon, \ldots, \varepsilon)} \!\!\!\!\!\!\!\!\!\! c_4 \cdots c_{m}.
\end{align*}
Now $\mathcal{T}(1,1,1,\varepsilon, \ldots, \varepsilon)=\{J\subset I:= \{3,\ldots,m\}\mid 3\in J\}$
and so Theorem~\ref{thm:6.1} yields
\begin{align*}
\dint\limits_{\M(1,1,1, \varepsilon, \ldots, \varepsilon)} \!\!\!\!\!\!\!\!\!\! c_4 \cdots c_{m} & =\sum_{J\in \mathcal{T}} (-1)^{\left(\sum_{i\in I\setminus J} 1 \right) + m - \vert J \vert } = \sum_{J\in \mathcal{T}} (-1)^{\vert I \vert - \vert J \vert + m - \vert J \vert } \\ & = \sum_{J\in \mathcal{T}} (-1)^{m-2+m}= \vert \mathcal{T} \vert = \sum_{i=0}^{m-3} \left( \begin{array}{cc} m-3 \\ i \end{array}\right)=2^{m-3}.
\end{align*}
For example, if $m=5$, $\M(\varepsilon,\varepsilon,1,1,1)$ is $S^2\times S^2$ and the Poincar\'e duals of $c_1$ and $c_2$ are $2[\{\text{pt}\}\times S^2]$ and $2[S^2\times \{\text{pt}\}]$. Hence, $\int_{\M(\alpha)} c_1c_2=4$ while $\int_{\M(\alpha)} c_1^2 = \int_{\M(\alpha)} c_2^2 = 0$. Moreover, $c_3=c_4=c_5=0$. 

In general, for $\M(\varepsilon,\ldots, \varepsilon,1,1,1)=\prod_{j=1}^{m-3} S^2_{\alpha_j}$, we have
$$
c_i = \text{PD}\left(2[\prod_{j\neq i}  S^2_{\alpha_j}]\right) \quad \text{for}\quad i=1,\ldots, m-3,
$$
while $c_{m-2}=c_{m-1}=c_m=0$.
\end{Example}

\section{Equilateral polygon spaces}\label{se:8}
Here we study the equilateral case corresponding to $\alpha_i=1$ for all $i$.   Since we want $\alpha$ to be generic we need $m$ to be odd. The intersection numbers for these spaces were computed by Kamiyama and Tezuka in \cite{KT}. Their result states the following.
\begin{Theorem}\label{thm:8.1}
Let $(d_1,\ldots,d_m)$ be a sequence of nonnegative integers with $\sum d_i=m-3$. Let $\beta_i, \varepsilon_i$ be such that $d_i=2\beta_i + \varepsilon_i$, where $\varepsilon_i=0$ or $1$. Then, defining
$$
\rho_{m,2k}:= (-1)^k \frac{\left( \begin{array}{c} \frac{m-3}{2} \\ k \end{array} \right) \left(\begin{array}{c} m-2 \\ \frac{m-1}{2} \end{array}\right)} { \left( \begin{array}{c} m-2  \\  2k + 1 \end{array} \right)},
$$
we have
\begin{enumerate}
\item if $\beta_i=0$ for $1\leq i \leq m$ then $\int c_1^{d_1}\cdots c_m^{d_m} = \rho_{m,0}$;
\item[]
\item if $\beta_i \neq 0$ for some $i$ then $\int c_1^{d_1}\cdots c_m^{d_m} = \rho_{m,2k}$, with $k=\beta_1+\cdots+ \beta_m$.
\end{enumerate}
\end{Theorem}
We will see how to obtain this result using Theorem~\ref{thm:6.1}.
\begin{proof}

{\bf 1.} First we see that if $\beta_i=0$ for $1\leq i\leq m$ then, since $\sum d_i=m-3$, we have
$$
\int_{\mathcal{M}_m}  c_1^{d_1}\cdots c_m^{d_m} = \int c_4 \cdots c_m.
$$
Moreover, 
$\mathcal{T}=\left\{ J\subset\{3,\ldots,m\}\mid\,\, \vert J \vert =\frac{m-1}{2}\right\}$ and so
\begin{align*}
\int_{\M_m} c_{4} \cdots c_m & = \sum_{J \in \mathcal{T}} (-1)^{\vert (I\setminus \{3\}) \setminus J\vert + m - \vert J  \vert} =  \sum_{J \in \mathcal{T}} (-1)^{\vert (I\setminus \{3\}) \setminus J\vert + \frac{m+1}{2} }\\ = & \sum_{J \in \mathcal{T} \text{s.t.} \, 3\in J} (-1)^{ \frac{m-3}{2} + \frac{m+1}{2} } +  \sum_{J \in \mathcal{T}  \text{s.t.}\, 3 \notin J}(-1)^{\frac{m-3}{2} - 1 + \frac{m+1}{2} } \\ = & \left( \begin{array}{c} m-3 \\ \frac{m-3}{2} \end{array} \right) (-1)^{m-1} +\left( \begin{array}{c} m-3 \\ \frac{m-1}{2} \end{array} \right) (-1)^{m-2} \\ = & (-1)^{m-1} \left\{  \left( \begin{array}{c} m-3 \\ \frac{m-3}{2} \end{array} \right) - \left( \begin{array}{c} m-3 \\ \frac{m-1}{2} \end{array} \right) \right\} = \frac{(m-3)!}{\left(\frac{m-1}{2}\right)!\left(\frac{m-3}{2}\right)!} = \rho_{m,0}.
\end{align*}

{\bf 2.} Let us now assume that $\beta_j\neq 0$ for some $j$. Since $d_1+\cdots+d_m=m-3$ we must have $d_j=0$ for some $j$. Then, using the fact that $c_1^2=c_2^2=\cdots=c_m^2$ (cf. \cite{HK2,B,K}) as well as  Remark~\ref{rk:3.1}, we have
$$
\int_{\mathcal{M}_m}  c_1^{d_1}\cdots c_m^{d_m} = \int_{\mathcal{M}_m} c_{2k+3}c_{2k+4}\cdots c_{m-1} c_m^{2k} 
$$
with $k=\beta_1+\cdots + \beta_m$. Now we have to consider two cases. 

If $2k=m-3$ then, since $\mathcal{T}=\left\{ J\subset I:=\{3,\ldots,m\}\mid\,\, \vert J \vert =\frac{m-1}{2}\right\}$, we have 
\begin{align}\label{eq:8.1}
\int_{\M_m} c_m^{2k} & = \sum_{J \in \mathcal{T}} (-1)^{\sum_{i\in (I\setminus J)} k_i + m - \vert J  \vert } =  (-1)^{\frac{m+1}{2}} \sum_{J \in \mathcal{T}} (-1)^{\sum_{i\in I\setminus J} k_i} \\\nonumber = &  (-1)^{\frac{m+1}{2}}  \vert \mathcal{T} \vert =  (-1)^{\frac{m+1}{2}} \left( \begin{array}{c} m-2 \\\nonumber \frac{m-1}{2} \end{array} \right) =  \rho_{m,m-3},
\end{align}
where we used the fact that all the exponents $k_i$ in \eqref{eq:8.1} are even ($k_i=0$ if $i\neq m$ and $k_m=2k$). 

If, on the other hand, $2k\neq m-3$ let us consider the set $A:=\{2k+3,\ldots,m-1\}\subset I$. Then
\begin{align*}
\int_{\M_m}c_{2k+3}c_{2k+4}\cdots c_{m-1} c_m^{2k}   & = \sum_{J \in \mathcal{T}} (-1)^{\sum_{i\in (I\setminus J)} k_i + m - \vert J  \vert } =  (-1)^{\frac{m+1}{2}} \sum_{J \in \mathcal{T}} (-1)^{\vert (I\setminus J) \cap A \vert}.
\end{align*}
Note that $\vert I \setminus J\vert =\frac{m-3}{2}$, $\vert I \setminus A\vert =2k+1$, $\vert A \vert = m-2k-3$ and that if there exist $j$ elements in $(I\setminus J) \cap (I\setminus A)$  then there exist $\frac{m-3}{2}-j$ elements in $(I \setminus J) \cap A$. Hence, 
\begin{align*}
\int_{\M_m}c_{2k+3}c_{2k+4}\cdots c_{m-1} c_m^{2k}   & = (-1)^{\frac{m+1}{2}} \sum_{j=0}^{2k+1}  \left( \begin{array}{c} 2k+1 \\ j \end{array} \right) \left( \begin{array}{c} m-2k-3 \\ \frac{m-3}{2}-j \end{array} \right) (-1)^{ \frac{m-3}{2}-j} \\ & =  \sum_{j=0}^{2k+1} (-1)^j  \left( \begin{array}{c} 2k+1 \\ j \end{array} \right) \left( \begin{array}{c} m-2k-3 \\ \frac{m-3}{2}-j \end{array} \right). 
\end{align*}
The result now follows from the following combinatorial identity, 
\begin{equation}
\label{eq:8.2}
\sum_{j=0}^{2k+1} (-1)^{j+k} \frac{\left( \begin{array}{c} 2k+1 \\ j \end{array} \right) \left( \begin{array}{c} m-2k-3 \\ \frac{m-3}{2}-j \end{array} \right) \left( \begin{array}{c} m-2 \\ 2k+1 \end{array} \right)}{ \left( \begin{array}{c} m-2 \\ \frac{m-1}{2} \end{array} \right)} = \left( \begin{array}{c} \frac{m-3}{2} \\ k \end{array} \right), 
\end{equation}
since then 
\begin{align*}
\int_{\M_m}c_{2k+3}c_{2k+4}\cdots c_{m-1} c_m^{2k}  =(-1)^k \frac{\left(\begin{array}{c} \frac{m-3}{2} \\ k \end{array}\right)\left(\begin{array}{c} m-2 \\ \frac{m-1}{2} \end{array}\right)}{\left(\begin{array}{c} m-2 \\ 2k+1 \end{array}\right)}=\rho_{2,2k}.
\end{align*}
Hence to finish the proof we just have to show  \eqref{eq:8.2}. For that we see that
\begin{align*}
& \sum_{j=0}^{2k+1}(-1)^{j+k}   \frac{ \left( \begin{array}{c} 2k+1 \\ j \end{array} \right) \left( \begin{array}{c} m-2k-3 \\ \frac{m-3}{2} - j \end{array} \right) \left( \begin{array}{c} m-2 \\ 2k+1 \end{array} \right) }{ \left( \begin{array}{c} m-2 \\ \frac{m-1}{2} \end{array} \right) } \\ & = \sum_{j=0}^{2k+1}(-1)^{j+k} \frac{\left( \frac{m-1}{2} \right)! \left(\frac{m-3}{2} \right) !}{j! (2k+1-j)!\left( \frac{m-3}{2} - j \right)! \left(\frac{m-3}{2} - 2k + j \right)!} \\ & =  \sum_{j=0}^{2k+1}(-1)^{j+k} \left( \begin{array}{c} \frac{m-3}{2} \\ j \end{array} \right) \left( \begin{array}{c} \frac{m-1}{2} \\ 2k+1-j \end{array} \right) = \left( \begin{array}{c} \frac{m-3}{2} \\ k \end{array} \right),
\end{align*}
where the last equality was obtained from the combinatorial Lemma~\ref{le:8.1} stated below, using $a=\frac{m-3}{2}$ and $b=2k+1$.
\end{proof}
\begin{Lemma}\label{le:8.1}
For any odd integer $b$ we have
$$
\sum_{j=0}^b (-1)^j \left(\begin{array}{c} a \\ j \end{array}\right) \left(\begin{array}{c} a+1 \\ b-j \end{array}\right) = (-1)^{\frac{b-1}{2}} \left( \begin{array}{c} a \\ \frac{b-1}{2} \end{array}\right).
$$
\end{Lemma}
\begin{proof} (Lemma~\ref{le:8.1})
This binomial identity can be easily proved using a generating function (cf. \cite{GKP} for a definition of generating function in this context). Indeed, since 
\begin{equation}\label{eq:8.3}
(z-1)^a(z+1)^{a+1}=(z^2-1)^a(1+z)
\end{equation}
and
\begin{align*}
& (z-1)^a  = \sum_{j\geq 0} (-1)^{a+j} \left(\begin{array}{c} a \\ j \end{array} \right) z^j, \quad \quad (z+1)^a  = \sum_{j\geq 0} \left(\begin{array}{c} a \\ j \end{array} \right) z^j, \\
& (z^2-1)^a(1+z)  =\sum_{j \geq 0} (-1)^{a+j} \left(\begin{array}{c} a \\ j \end{array} \right) z^{2j}+ \sum_{j\geq 0} (-1)^{a+j} \left(\begin{array}{c} a \\ j \end{array} \right) z^{2j+1},
\end{align*}
the result follows from equating coefficients of $z^b$ in \eqref{eq:8.3} for an odd $b$.
\end{proof}
\begin{Example}
If $m=5$ we obtain 
$$
\int_{\M_5} c_i^2 = \rho_{5,2} = -3\,\,\, \text{and} \,\, \int_{\M_5} c_i c_j = \rho_{5,0}= 1, \, \text{if}\, i \neq j.
$$
Incidentally, note that these values are equal to the ones obtained in Example~\ref{ex:7.1}. This is not surprising since $\M_5$ is diffeomorphic to $\M(4,3,4,3,4)$. Indeed, by Remark~\ref{rk:2.1},  $\M_5$ is diffeomorphic to $\M(\alpha_\varepsilon)$ with $\alpha_\varepsilon=(1+\varepsilon, 1,1,1,1+ \varepsilon)$ for small values of $\varepsilon$, and $\M(\alpha_\varepsilon)$ is a toric manifold whose moment polytope has $7$ edges just like the one for $\M(4,3,4,3,4)$ depicted in Figure~\ref{fig:7.1} (the three spaces  are all diffeomorphic to $(S^2 \times S^2)\# 3 \overline{\C P^2}$). 
\end{Example}
The space of equilateral polygons $\M_m$ admits a natural action of the symmetric group $Sym_m$. Moreover, the quotient $\M_m/Sym_m$ can be seen as a compactification of the moduli space of $m$ unordered points in $\C P^1$ as well as the space of $m$-times punctured genus zero algebraic curves. The cohomology ring of this space was computed by Brion \cite{B}, Klyachko \cite{K} and Hausmann and Knutson \cite{HK2}. Since 
$$
H^*(\M_m/Sym_m;\Q) \cong H^*(\M_m; \Q)^{Sym_m},
$$
can be identified with the invariant part of $H^*(\M_m;\Q)$, and  $H^*(\M_m;\Q)$ is generated in degree $2$, we just have to study the action of $Sym_m$ on $H^2(\M_m)$. Hausmann and Knutson \cite{HK2} prove that this $Sym_m$-invariant part of   $H^*(\M_m; \Q)$ is generated by $\sigma_1$, the first invariant symmetric polynomial in the classes $c_i$, and by $c_m^2$ (or any other $c_i^2$ since they are all equal to the Pontrjagin class of the principal $SO(3)$-bundle $A_m \to \M_m$ with $A_m:=\{(\vec{v}_1,\ldots,\vec{v}_m)\in (\R^3)^m\mid \sum_{i=1}^m \vec{v}_i = 0,\,\, \vert \vec{v}_1\vert= \cdots = \vert \vec{v}_m\vert\}$), a generator of degree $4$, with no relations up to degree $m-3$. The intersection numbers $\sigma_1^k \cdot c_m^{m-3-k}$ for an even $k$ can be obtained from Theorem~\ref{thm:8.1}.
First we see that
$$
\sigma_1^k = \sum_{k_1+\cdots+k_m=k} \frac{k!}{k_1!\cdots k_m!} c_1^{k_1}\cdots c_m^{k_m}.
$$
Then
\begin{align*}
\int_{\M_m/Sym_m} \sigma_1^k \cdot c_m^{m-3-k} & = \sum_{k_1+\cdots+k_m=k} \frac{k!}{k_1!\cdots k_m!} \int_{\M_m}  c_1^{k_1}\cdots c_{m-1}^{k_{m-1}}c_m^{k_m+m-3-k} \\ & = \sum_{k_1+\cdots+k_m=k} \frac{k!}{k_1!\cdots k_m!} \rho_{m,2(\beta_1+\cdots+\beta_{m-1}+\tilde{\beta}_m)},
\end{align*}
where the $\beta_i$ are such that $k_i=2\beta_i + \varepsilon_i$, $\varepsilon_i=0$ or $1$ and $\tilde{\beta}_m=\beta_m + \frac{m-3-k}{2}$. Hence, 
$$
2(\beta_1+\cdots+\beta_{m-1}+\tilde{\beta}_m)=m-3- \text{number of odd $k_i$'s}
$$
and so, since the number of odd $k_i$'s must be even ($k$ is even), we have
$$
\int_{\M_m/Sym_m} \sigma_1^k \cdot c_m^{m-3-k} = \sum_{j=0}^{\frac{k}{2}} \sum_{ \small{\begin{array}{c} k_1 + \cdots + k_m =k \,\, \text{s.t.} \\ 2j\, \text{of the $k_i$'s are odd}\end{array}}} \frac{k!}{k_1!\cdots k_m!} \, \rho_{m,m-3-2j}
$$
and we conclude the following result.
\begin{Proposition} For an even integer $k$
\begin{align*}
& \int_{\M_m/Sym_m}  \sigma_1^k \cdot c_m^{m-3-k}  = \\ &=  (-1)^{\frac{m-3}{2}}  \left( \begin{array}{c} m-2 \\ \frac{m-1}{2} \end{array}\right)  \sum_{j=0}^{\frac{k}{2}} (-1)^j \frac{ \left( \begin{array}{c} \frac{m-3}{2} \\ j  \end{array} \right)}{ \left( \begin{array}{c} m-2 \\ 2j \end{array}\right)} \sum_{ \small{\begin{array}{c} k_1+\cdots+k_m=k \, \text{s.t.} \\ 2j\, \text{of the $k_i$'s are odd}\end{array}}} \left(\begin{array}{c} k \\ k_1, \ldots, k_m \end{array}\right).
\end{align*}
\end{Proposition}

\section{Moduli space of flat connections} 
\label{se:9}

Polygon spaces  can be identified with moduli spaces of flat $SU(2)$-connections on the $m$-punctured sphere. In this section we compare the existing formulas for cohomology intersection pairings in the context of moduli spaces of flat connections with our explicit  expression of Theorem~\ref{thm:6.1}.

Let us consider closed polygons in $S^3$ with vertices $v_1,\ldots,v_n$ joined by edges $e_1,\ldots,e_n$ where each $e_i$ is the geodesic arc from $v_i$ to $v_{i+1}$. The length $\alpha_i$ of an edge is then the length of this geodesic arc. Let us denote by $\M(\alpha)^{S^3}$ the moduli space of closed polygons in $S^3$ with side-lengths $\alpha$ modulo orientation-preserving isometries, that is modulo $SO(4)$. This space can be identified with the moduli space of flat $SU(2)$ connections on a punctured sphere with fixed holonomies around the punctures. Indeed, let us denote by $\M(S_m,\alpha)$ the moduli space of flat $SU(2)$ connections on the $m$-punctured sphere $S_m:=S^2 \setminus \{p_1, \ldots, p_m\}$ modulo gauge equivalence such that the holonomy around $p_i$ is conjugate to 
$$
A_{\alpha_i}:=\left( \begin{array}{cc} e^{i \pi \alpha_i} & 0 \\ 0 & e^{-i \pi \alpha_i} \end{array}\right).
$$
Then,
$$
\M(S_m,\alpha) \cong \mathcal{R}(m,\alpha)/SU(2),
$$
where 
$$
\mathcal{R}(m,\alpha)= \{(g_1,\ldots, g_m)\in SU(2)^m\mid g_1\ldots g_m = Id,\, \text{tr} g_i = 2 \cos \pi \alpha_i, \, i=1,\ldots,m\}
$$
is the space of representations from $\pi_1(S_m)$ to $SU(2)$ such that the image of the loop around the puncture $p_i$ is conjugate to $A_{\alpha_i}$ and $SU(2)$ acts diagonally.

Moreover, $\mathcal{R}(m,\alpha)$ can be identified with the space $\mathcal{P}(\alpha)^{S^3,0}$ of based polygons in $S^3$ that have the first vertex fixed $v_1=Id \in SU(2)\cong S^3$. Each element $(g_1,\ldots,g_m)\in \mathcal{R}(m,\alpha)$ is identified with the based polygon that has vertices
$$
v_1=Id, \quad v_2=g_1, \quad v_3=g_1g_2, \quad \cdots \quad v_i=g_1g_2\cdots g_{i-1}, \quad i=1,\ldots, m.
$$
Then, 
$$
\M(S_m,\alpha) \cong \mathcal{R}(m,\alpha)/SU(2)\cong \mathcal{P}(\alpha)^{S^3,0}/SO(3) \cong \M(\alpha)^{S^3}/SO(4).
$$
A formula for the symplectic volume of the space  $\M(S_m,\alpha)$ was first obtained by Witten in \cite{Wi2} and proved rigorously by Jeffrey and Weitsman in \cite{JW}. It states the following.
\begin{Theorem}(Witten, Jeffrey-Weitsman) 
\begin{equation}\label{eq:9.1}
\text{Vol} (\M(S_m,\alpha))= \frac{4}{\pi^{m-2}} \sum_{k=1}^\infty \frac{\prod_{i=1}^m \sin (k \pi \alpha_i)}{k^{m-2}}.
\end{equation}
\end{Theorem} 
Later in \cite{V} Vu The Khoi obtained a closed form expression for this volume in terms of Bernoulli polynomials.
Noting that the moduli spaces of polygons in $S^3$ and in the Euclidean space $\R^3$ with the same side lengths are symplectomorphic provided that the $\alpha_i$ are sufficiently small (cf. Theorem 6.6. in \cite{J} and \cite{KM}), and that multiplying $\alpha$ by a scalar $\lambda>0$ we get $\text{vol}( \M(\lambda \alpha))=\lambda^{m-3} \text{vol} (\M(\alpha))$, he deduces an expression for the volume of $\M(\alpha)$.
\begin{Proposition}{(Vu The Khoi)}\label{prop:9.1}
The symplectic volume of the moduli space $\M(\alpha)$ for $\alpha=(\alpha_1,\ldots, \alpha_m)$ is given by
\begin{equation}
\nonumber
\text{Vol}(\M (\alpha)) = - \frac{1}{4 (m-3)!} \sum_{R \subset \{1,\ldots, m \}} (-1)^{\vert R \vert } \sgn (\sum_{i \in R} \alpha_i - \sum_{i \notin R} \alpha_i)(\sum_{i \in R} \alpha_i - \sum_{i \notin R} \alpha_i)^{m-3} 
\end{equation}
if $m$ is even and by
\begin{equation}
\nonumber
\text{Vol}(\M(\alpha))= - \frac{1}{2 (m-3)!} \sum_{R \subset \{1,\ldots,m\} \, \text{s.t.} \, \vert R \vert \, \text{odd}}  \sgn (\sum_{i \in R} \alpha_i - \sum_{i \notin R} \alpha_i)(\sum_{i \in R} \alpha_i - \sum_{i \notin R} \alpha_i)^{m-3} 
\end{equation}
if $m$ is odd.
\end{Proposition}
A similar expression was obtained independently by Mandini in \cite{M} using localization theorems in equivariant cohomology and an equivariant integration formula for symplectic quotients by non-abelian groups.

Using the Witten-Jeffrey-Weitsman expression for the volume of $\M (S_m,\alpha)$ Yoshida obtains in \cite{Y} a generating function for cohomology intersection pairings of the moduli space of flat connections.
\begin{Theorem}(Yoshida)
Let $t_1, \ldots, t_m \in \R$. Then,
\begin{align*}
\sum_{k_1,\ldots, k_m \geq 0} \frac{t_1^{k_1}}{k_1!} \cdots \frac{t_m^{k_1}}{k_m!} \int_{\M(S_m, \alpha)} c_1^{k_1} \cdots c_1^{k_m} = \frac{4}{\pi^{m-2}} \sum_{l=1}^\infty \frac{ \prod_{j=1}^m \sin (\pi \, l\, (\alpha_j + x_j))}{l^{m-2}}
\end{align*}
where $\sum_{i=1}^m k_i=m-3$.
\end{Theorem}
The cohomology intersection pairings $\int_{\M(S_m, \alpha)} c_1^{k_1} \cdots c_1^{k_m}$ can then be obtained from the above expression \eqref{eq:9.1} for the volume by taking the appropriate derivatives
\begin{equation}
\nonumber
\int_{\M(S_m, \alpha)} c_1^{k_1} \cdots c_1^{k_m}= \frac{\partial^{k_1}}{\partial t_1^{k_1}} \cdots \frac{\partial^{k_m}}{\partial t_m^{k_m}} \text{Vol}(\M(S_m, \alpha + t))\vert_{t=0}.
\end{equation}
Using this formula and the identification of $\M(S_m,\alpha)$ with $\M(\alpha)$ for small values of $\alpha$ one obtains a formula for intersection pairings in $\M(\alpha)$ by taking derivatives of the  expressions in Proposition~\ref{prop:9.1}.
\begin{Theorem}\label{thm:9.3}
\begin{equation}
\label{eq:9.2}
\int_{\M(S_m, \alpha)} c_1^{k_1} \cdots c_1^{k_m} = \frac{1}{4} 
\sum_{R \subset \{1,\ldots, m \}} (-1)^{\vert R \vert + 1 + \sum_{i \notin R} k_i} \sgn (\sum_{i \in R} \alpha_i - \sum_{i \notin R} \alpha_i)
\end{equation}
if $m$ is even and 
\begin{equation}
\label{eq:9.3}
\int_{\M(S_m, \alpha)} c_1^{k_1} \cdots c_1^{k_m} = \frac{1}{2} \sum_{R \subset \{1,\ldots,m\} \, \text{s.t.} \, \vert R \vert \, \text{odd}} (-1)^{1+ \sum_{i \notin R} k_i} \sgn (\sum_{i \in R} \alpha_i - \sum_{i \notin R} \alpha_i)
\end{equation}
if $m$ is odd.
\end{Theorem}
We will now see how this formula is equivalent to our explicit expression of Theorem~\ref{thm:6.1}.
\begin{proof} Let us assume without loss of generality that $\alpha_1 > \alpha_2$ and let $I=\{3,\ldots,m\}$. Notice that
\begin{align*}
\{ & R \subset \{1, \ldots, m\} \}  = \{ R \subset I \mid  R\in \mathcal{T}(\alpha) \} \bigcup \{ R \subset I \mid  R\notin \mathcal{T}(\alpha) \}  \bigcup \\ & \bigcup  \{R=\{1,2\}\cup J\mid J\subset I \}  \bigcup \{R= \{1\} \cup J\mid  J\subset I \}  \bigcup \{R= \{2\} \cup J\mid J\subset I \}.
\end{align*}
Let us denote by $S_R$ the difference
$$
S_R:= \sum_{i \in R} \alpha_i - \sum_{i \notin R} \alpha_i
$$
and by $l_J$ the sum 
$$
l_J:=\sum_{i=3}^m (-1)^{\chi_{I \setminus J}(i)}\alpha_i,
$$
whenever $J\subset I$.

If $R\subset I$ and $R\in \mathcal{T}(\alpha)$ then $S_R= l_R - \alpha_1-\alpha_2 \leq 0$
since for triangular sets $l_R \leq \alpha_1 + \alpha_2$. If $R\subset I$ is not triangular then, if $l_R \leq 0$ or $0< l_R < \alpha_1 - \alpha_2$, we have  $S_R < 0$. If, however, $l_R> \alpha_1 + \alpha_2$ then $S_R>0$.

If $R=\{1,2\} \cup J$ with $J\subset I$  we have $\vert R \vert = \vert J \vert + 2$ and $S_R=l_J + \alpha_1 +\alpha_2$. If in addition $J\in \mathcal{T}(\alpha)$ then $S_R>0$. If, on the other hand, $J\notin \mathcal{T}(\alpha)$ then $S_R>0$ whenever $0 > l_J > -(\alpha_1 + \alpha_2)$ or $l_J>0$ and negative otherwise.

If $R=\{1\} \cup J$ with $J\subset I$ then $\vert R \vert = \vert J \vert + 1$ and $S_R=l_J + \alpha_1 - \alpha_2$. If $J\in \mathcal{T}(\alpha)$ then $S_R>0$. If, however, $J\notin \mathcal{T}(\alpha)$ then $S_R>0$ whenever $0 > l_J > -(\alpha_1 - \alpha_2)$ or $l_J>0$ and negative otherwise.

Similarly, if $R=\{2\} \cup J$ with $J\subset I$  we have $\vert R \vert = \vert J \vert + 1$ and $S_R=l_J + \alpha_2 - \alpha_1$. If in addition $J\in \mathcal{T}(\alpha)$ then $S_R>0$ (since by the triangular inequalities $l_J\geq \vert \alpha_2 - \alpha_1\vert$). If, on the other hand, $J\notin \mathcal{T}(\alpha)$ then $S_R>0$ if $l_J > \alpha_1 + \alpha_2$  and negative otherwise.

Let us assume  first that $m$ is even.
Putting the above information  together, the RHS of  \eqref{eq:9.2} is equal to
\begin{align*}
& \frac{1}{4} \left( \sum_{J \in \mathcal{T}(\alpha)} (-1)^{\vert J \vert +  \sum_{i \in I\setminus J} k_i}\left((-1)^{ k_1 + k_2} -1 + (-1)^{k_2} + (-1)^{k_1}\right) + \right. \\ 
& +  \sum_{J \notin \mathcal{T}(\alpha)\,\,\text{s.t.}\,\, l_J> \alpha_1 + \alpha_2} (-1)^{\vert J \vert +  \sum_{i \in I\setminus J} k_i}\left((-1)^{1+k_1+k_2}-1+(-1)^{k_2}+(-1)^{k_1}\right) +  
\\ 
& + \sum_{J \notin \mathcal{T}(\alpha)\,\,\text{s.t.}\,\, 0<l_J< \alpha_1 - \alpha_2} (-1)^{\vert J \vert +  \sum_{i \in I\setminus J} k_i}\left((-1)^{k_1+k_2}-1+(-1)^{k_2}-(-1)^{k_1}\right) + \\ 
& + \sum_{J \notin \mathcal{T}(\alpha)\,\,\text{s.t.}\,\, l_J<-(\alpha_1 + \alpha_2)} (-1)^{\vert J \vert +  \sum_{i \in I\setminus J} k_i }\left((-1)^{k_1+k_2}+1-(-1)^{k_2}-(-1)^{k_1}\right) + \\ 
 & + \sum_{J \notin \mathcal{T}(\alpha)\,\,\text{s.t.}\,\, -(\alpha_1 + \alpha_2)<l_J<\alpha_2 - \alpha_1} \hspace{-1.1cm} (-1)^{\vert J \vert +  \sum_{i \in I\setminus J} k_i }\left((-1)^{k_1+k_2}-1-(-1)^{k_2}-(-1)^{k_1}\right) + \\        
 & \left. + \sum_{J \notin \mathcal{T}(\alpha)\,\,\text{s.t.}\,\, \alpha_2 - \alpha_1 <l_J < 0} (-1)^{\vert J \vert +  \sum_{i \in I\setminus J} k_i }\left((-1)^{k_1+k_2}-1+(-1)^{k_2}-1\right)\right). 
\end{align*} 
Assuming as in Theorem~\ref{thm:6.1}  that $k_{m-l},\ldots, k_m \in \Z_+$  for  some integer $l$ and $k_1=\cdots = k_{m-l-1}=0$ then  $k_{m-l}+\cdots + k_m = m-3$ and we must have $m-l \geq 3$. Thus  $k_1=k_2=0$. Hence, the RHS of  \eqref{eq:9.2} is equal to 
\begin{equation}
\label{eq:9.4}
\frac{1}{2} \sum_{J \in \mathcal{T}(\alpha)} (-1)^{\vert J \vert +  \sum_{i \in I\setminus J} k_i} - \frac{1}{2}\hspace{-1cm} \sum_{\tiny{\begin{array}{c} J \notin \mathcal{T}(\alpha)\,\,\text{s.t.}\\  -(\alpha_1 + \alpha_2)<l_J<\alpha_2 - \alpha_1 < 0\end{array}}}\hspace{-1cm} (-1)^{\vert J \vert +  \sum_{i \in I\setminus J} k_i}. 
\end{equation} 
If $J\subset I$ is not triangular and 
\begin{equation}
\label{eq:9.5} -(\alpha_1 + \alpha_2 ) < l_J < \alpha_2 -\alpha_1 < 0
\end{equation} 
then, the complement $J^\prime:= I \setminus J$ satisfies  $\alpha_1 -\alpha_2 < l_{J^\prime} < \alpha_1 + \alpha_2$ (since $l_{J^\prime}=-l_J$), implying that $J^\prime$ is triangular. Conversely, if a subset of $I$ is triangular its complement in $I$ satisfies \eqref{eq:9.5}. Hence, since  $\vert J\vert = m-2 -\vert J^\prime \vert$, the sum  \eqref{eq:9.4} becomes equal to    
\begin{align*}
& \frac{1}{2} \sum_{J \in \mathcal{T}(\alpha)} (-1)^{\vert J \vert +  \sum_{i \in I\setminus J} k_i} - \frac{1}{2} \sum_{J^\prime \in \mathcal{T}(\alpha)} (-1)^{m- \vert J^\prime  \vert +  (m-3 - \sum_{i \in I\setminus J^\prime} k_i)} \\ &  = \sum_{J \in  \mathcal{T}(\alpha)} (-1)^{m- \vert J \vert +  \sum_{i \in I\setminus J} k_i},  
\end{align*} 
which is our formula of Theorem~\ref{thm:6.1}.

Similarly, if $m$ is odd, the RHS of \eqref{eq:9.3} is equal to 
\begin{align*}
&  \frac{1}{2}  \left( \sum_{\tiny{ J \in \mathcal{T}(\alpha), \,\, \vert J  \vert \,\,\text{odd}}} \right.  (-1)^{ \sum_{i \in I\setminus J} k_i}\left((-1)^{ k_1 + k_2} -1\right) \, \, + \\ & + \hspace{-.5cm} \sum_{\tiny{\begin{array}{c} J \notin \mathcal{T}(\alpha),\,\, \vert J  \vert \,\,\text{odd} \\ l_J > \alpha_1 + \alpha_2 \end{array}}}  \hspace{-.8cm}(-1)^{\sum_{i \in I\setminus J} k_i}\left((-1)^{1+k_1+k_2}-1\right) +\hspace{-.5cm}  \sum_{\tiny{\begin{array}{c} J \notin \mathcal{T}(\alpha),\,\, \vert J  \vert \,\,\text{odd} \\  0<l_J < \alpha_1 - \alpha_2\end{array}}}  \hspace{-.8cm} (-1)^{\sum_{i \in I\setminus J} k_i}\left((-1)^{k_1+k_2}-1\right) + \\  & +
 \hspace{-.5cm}  \sum_{\tiny{\begin{array}{c}  J \notin \mathcal{T}(\alpha),\,\, \vert J  \vert \,\,\text{odd} \\  l_J < -(\alpha_1 + \alpha_2)\end{array}}}  \hspace{-.8cm} (-1)^{\sum_{i \in I\setminus J} k_i }\left((-1)^{k_1+k_2}+1\right) +  \hspace{-.5cm}  \sum_{\tiny{\begin{array}{c} J \notin \mathcal{T}(\alpha),\,\, \vert J  \vert \,\,\text{odd} \\  -(\alpha_1 + \alpha_2)<l_J <0 \end{array}}} \hspace{-.8cm} (-1)^{ \sum_{i \in I\setminus J} k_i }\left((-1)^{k_1+k_2}-1\right) + \\  & +
\sum_{\tiny{ J \in \mathcal{T}(\alpha), \,\, \vert J  \vert \,\,\text{even}}} (-1)^{ \sum_{i \in I\setminus J} k_i}\left((-1)^{1 + k_2} + (-1)^{1 + k_1}\right) +  \\ &  + \hspace{-.5cm} \sum_{\tiny{\begin{array}{c} J \notin \mathcal{T}(\alpha),\,\, \vert J  \vert \,\,\text{even} \\ l_J > \alpha_1 + \alpha_2 \end{array}}}  \hspace{-1cm}(-1)^{\sum_{i \in I\setminus J} k_i}\left((-1)^{1+k_2}+(-1)^{1+k_1}\right) +\hspace{-1cm}  \sum_{\tiny{\begin{array}{c} J \notin \mathcal{T}(\alpha),\,\, \vert J  \vert \,\,\text{even} \\  0 < l_J < \alpha_1 - \alpha_2\end{array}}}  \hspace{-1cm} (-1)^{\sum_{i \in I\setminus J} k_i}\left((-1)^{1+k_2}+(-1)^{k_1} \right) + \\  & \left. +
 \hspace{-.5cm}  \sum_{\tiny{\begin{array}{c}  J \notin \mathcal{T}(\alpha),\,\, \vert J  \vert \,\,\text{even} \\  l_J < \alpha_2 - \alpha_1<0)\end{array}}}  \hspace{-.8cm} (-1)^{\sum_{i \in I\setminus J} k_i }\left((-1)^{k_2}+(-1)^{k_1}\right) +  \hspace{-1cm}  \sum_{\tiny{\begin{array}{c} J \notin \mathcal{T}(\alpha),\,\, \vert J  \vert \,\,\text{odd} \\  \alpha_2 - \alpha_1 < l_J <0 \end{array}}} \hspace{-1cm} (-1)^{ \sum_{i \in I\setminus J} k_i }\left((-1)^{1+k_2}+(-1)^{k_1}\right) \right) = \\
& =- \hspace{-.5cm} \sum_{\tiny{\begin{array}{c} J \notin \mathcal{T}(\alpha),\,\, \vert J  \vert \,\,\text{odd} \\ l_J> \alpha_1 + \alpha_2 \end{array}}}  \hspace{-.8cm}(-1)^{\sum_{i \in I\setminus J} k_i} + \hspace{-.5cm}  \sum_{\tiny{\begin{array}{c}  J \notin \mathcal{T}(\alpha),\,\, \vert J  \vert \,\,\text{odd} \\ l_J <-(\alpha_1 + \alpha_2)\end{array}}} \hspace{-.8cm} (-1)^{\sum_{i \in I\setminus J} k_i }  - \hspace{-.5cm} \sum_{\tiny{J \in \mathcal{T}(\alpha),\,\, \vert J  \vert \,\,\text{even}}}  \hspace{-.8cm}(-1)^{\sum_{i \in I\setminus J} k_i} - \\ & -  \hspace{-.5cm} \sum_{\tiny{\begin{array}{c}  J \notin \mathcal{T}(\alpha),\,\, \vert J  \vert \,\,\text{even} \\ l_J > \alpha_1 + \alpha_2 \end{array}}} \hspace{-.8cm} (-1)^{\sum_{i \in I\setminus J} k_i } + \hspace{-.5cm}  \sum_{\tiny{\begin{array}{c}  J \notin \mathcal{T}(\alpha),\,\, \vert J  \vert \,\,\text{even} \\ -(\alpha_1 + \alpha_2) < l_J < \alpha_2 - \alpha_1 \end{array}}} \hspace{-.8cm} (-1)^{\sum_{i \in I\setminus J} k_i } + \hspace{-.5cm}  \sum_{\tiny{\begin{array}{c}  J \notin \mathcal{T}(\alpha),\,\, \vert J  \vert \,\,\text{even} \\ l_J < -(\alpha_1 + \alpha_2) \end{array}}} \hspace{-.8cm} (-1)^{\sum_{i \in I\setminus J} k_i }.
\end{align*} 
If $J\subset I$ is not triangular, $\vert J \vert$ is even and $l_J < -(\alpha_1 + \alpha_2)$, then the complement $J^\prime:= I \setminus J$  satisfies $l_{J^\prime} > \alpha_1 + \alpha_2$ and  $\vert J^\prime \vert$ is odd and vice-versa. Moreover, if $J\subset I$ is not triangular, $\vert J \vert$ is even and  $-(\alpha_1 + \alpha_2)< l_J <\alpha_2-\alpha_1 <0$, then $J^\prime$ is triangular, $\vert J^\prime \vert$ is odd  and vice-versa.  Hence, the RHS of \eqref{eq:9.3} is equal to
\begin{align*}
 - \sum_{\tiny{\begin{array}{c} J \in \mathcal{T}(\alpha)\\  \vert J  \vert \,\,\text{even}\end{array}}}  (-1)^{\sum_{i \in I\setminus J} k_i}  + \sum_{\tiny{\begin{array}{c} J \in \mathcal{T}(\alpha) \\  \vert J  \vert \,\,\text{odd}\end{array}}} (-1)^{\sum_{i \in I\setminus J} k_i} = \sum_{\tiny{ J \in \mathcal{T}(\alpha)}} (-1)^{m- \vert J \vert + \sum_{i \in I\setminus J} k_i}
\end{align*} 
which again is our formula from Theorem~\ref{thm:6.1}.
\end{proof}

\section{Other explicit formulas for intersection pairings} 
\label{se:10}
Using the commutativity of geometric quantization and symplectic reduction Takakura obtains in \cite{T} an explicit formula for cohomology intersection pairings on arbitrary polygon spaces. Later  Konno obtains an equivalent expression (see  \cite{Konno}), although written in a different basis, using algebro-geometric methods. Their result is the following.
\begin{Theorem}(Konno, Takakura) 

Let $\alpha=(\alpha_1,\ldots, \alpha_m)$ be generic and let
$(k_1,\ldots,k_m)$ be a sequence of nonnegative integers with $\sum_{i=1}^m k_i=m-3$. Let $\mathcal{S}(\alpha)$ be the family of sets  $R \subset \{1, \ldots, m \}$ for which $\sum_{i \in R} \alpha_i - \sum_{i \notin R} \alpha_i < 0$. Then we have
\begin{equation}
\label{eq:10.1}
\int_{\M(\alpha)} c_1^{k_1} \cdots c_m^{k_m} = -\frac{1}{2} \sum_{R \in \mathcal{S}(\alpha)} (-1)^{\vert R \vert + \sum_{i \in R} k_i}.
\end{equation}
\end{Theorem}
We will now see how this formula is equivalent to Theorem~\ref{thm:6.1} although our formula uses a smaller  family of  sets (the triangular sets $\mathcal{T}(\alpha)$) which is always contained in $\mathcal{S}(\alpha)$.
\begin{proof} Let us again assume without loss of generality that $\alpha_1 > \alpha_2$ and let $I=\{3,\ldots,m\}$. Notice that
\begin{align*}
& \mathcal{S}(\alpha)  = \{ R \subset I \mid  R\in \mathcal{T}(\alpha) \cap \mathcal{S}(\alpha) \} \bigcup \{ R \subset I \mid  R\notin \mathcal{T}(\alpha)\, \text{and}\, R \in \mathcal{S}(\alpha) \}  \bigcup \\ & \bigcup  \{R=\{1,2\}\cup J\mid J\subset I \, \text{and}\, R \in \mathcal{S}(\alpha)\}  \bigcup \{R= \{1\} \cup J\mid  J\subset I \, \text{and}\, R \in \mathcal{S}(\alpha) \} \bigcup   \\ & \bigcup \{R= \{2\} \cup J\mid J\subset I\, \text{and}\, R \in \mathcal{S}(\alpha) \}.
\end{align*}
Again let $S_R$ and $l_J$ be respectively equal to 
$$
S_R:= \sum_{i \in R} \alpha_i - \sum_{i \notin R} \alpha_i \,\,\,\, \text{and} \,\,\,\, l_J:=\sum_{i=3}^m (-1)^{\chi_{I \setminus J}(i)}\alpha_i,
$$
whenever $J\subset I$. Note that the elements of $\mathcal{S}(\alpha)$ are precisely the subsets of $\{1, \ldots, m\}$ for which $S_R<0$.

As we saw before in the proof of Theorem~\ref{thm:9.3}, if $R\subset I$ and $R\in \mathcal{T}(\alpha)$ then $S_R\leq 0$ and so $R\in \mathcal{S}(\alpha)$. On the other hand, if $R\subset I$ is not triangular then $R\in \mathcal{S}(\alpha)$ if and only if $l_R \leq 0$ or $0< l_R < \alpha_1 - \alpha_2$. 

If $R=\{1,2\} \cup J$ with $J\subset I$,  we have $\vert R \vert = \vert J \vert + 2$ and $S_R=l_J + \alpha_1 +\alpha_2$, and then $R\in \mathcal{S}(\alpha)$ if and only if  $J\notin \mathcal{T}(\alpha)$ and  $l_J < -(\alpha_1 + \alpha_2)$.

If $R=\{1\} \cup J$ with $J\subset I$ then $\vert R \vert = \vert J \vert + 1$ and $S_R=l_J + \alpha_1 - \alpha_2$ and so $R\in \mathcal{S}(\alpha)$ exactly when  $J\notin \mathcal{T}(\alpha)$ and  $l_J < -(\alpha_1 - \alpha_2)$.

Similarly, if $R=\{2\} \cup J$ with $J\subset I$,  we have $\vert R \vert = \vert J \vert + 1$ and $S_R=l_J + \alpha_2 - \alpha_1$, and so $R\in \mathcal{S}(\alpha)$ if and only if $J\notin \mathcal{T}(\alpha)$ and  $l_J < \alpha_1 + \alpha_2$.

Putting this information  together, the RHS of  \eqref{eq:10.1} is equal to
\begin{align*}
& - \frac{1}{2} \left( \sum_{J \in \mathcal{T}(\alpha)} (-1)^{\vert J \vert +  \sum_{i \in J} k_i} 
+ \hspace{-.5cm} \sum_{J \notin \mathcal{T}(\alpha)\,\,\text{s.t.}\,\, l_J < \alpha_1 - \alpha_2} \hspace{-1.2cm}(-1)^{\vert J \vert +  \sum_{i \in J} k_i} 
+ \hspace{-.5cm} \sum_{J \notin \mathcal{T}(\alpha)\,\,\text{s.t.}\,\, l_J< -(\alpha_1 + \alpha_2)}  \hspace{-1.3cm}(-1)^{\vert J \vert +  \sum_{i \in J} k_i} \right. \\ &\left. 
-\hspace{-.5cm} \sum_{J \notin \mathcal{T}(\alpha)\,\,\text{s.t.}\,\, l_J<-(\alpha_1 - \alpha_2)}\hspace{-1.2cm} (-1)^{\vert J \vert +  \sum_{i \in J} k_i} \hspace{.3cm}
- \hspace{-.5cm}\sum_{J \notin \mathcal{T}(\alpha)\,\,\text{s.t.}\,\, l_J< \alpha_1 + \alpha_2} \hspace{-1.1cm} (-1)^{\vert J \vert +  \sum_{i \in  J} k_i }\right) = \\
& - \frac{1}{2} \left( \sum_{J \in \mathcal{T}(\alpha)} (-1)^{\vert J \vert +  \sum_{i \in J} k_i} \hspace{.5cm}-\hspace{-.5cm}  \sum_{\tiny{\begin{array}{c} J \notin \mathcal{T}(\alpha)\,\,\text{s.t.}\\  -(\alpha_1 + \alpha_2)<l_J<\alpha_2 - \alpha_1 < 0\end{array}}}\hspace{-1.2cm} (-1)^{\vert J \vert +  \sum_{i \in J} k_i}\right).
\end{align*} 
As we saw before, if $J\subset I$ is not triangular and 
\begin{equation}
\label{eq:10.2} -(\alpha_1 + \alpha_2 ) < l_J < \alpha_2 -\alpha_1 < 0,
\end{equation} 
then the complement $J^\prime:= I \setminus J$ satisfies  $\alpha_1 -\alpha_2 < l_{J^\prime} < \alpha_1 + \alpha_2$, implying that $J^\prime$ is triangular. Conversely, if a subset of $I$ is triangular its complement in $I$ satisfies \eqref{eq:10.2}. Hence, since  $\vert J\vert = m-2 -\vert J^\prime \vert$ and  $\sum_{i \in I\setminus J} k_i = (m-3) - \sum_{i \in J} k_i$, the RHS of  \eqref{eq:10.1} becomes equal to    
\begin{align*}
& \frac{1}{2}\left( \sum_{J \in \mathcal{T}(\alpha)} (-1)^{\vert J \vert +  \sum_{i \in I\setminus J} k_i + m} + \sum_{J^\prime \in \mathcal{T}(\alpha)} (-1)^{m- \vert J^\prime  \vert +  \sum_{i \in I\setminus J^\prime} k_i} \right)\\ &  = \sum_{J \in  \mathcal{T}(\alpha)} (-1)^{m- \vert J \vert +  \sum_{i \in I\setminus J} k_i},  
\end{align*} 
which is our formula of Theorem~\ref{thm:6.1}.
\end{proof}

\end{document}